\documentstyle[12pt,psfig]{amsart}
\makeatletter
\def\eqalign#1{\null\,\vcenter{\openup\jot\m@th
  \ialign{\strut\hfil$\displaystyle{##}$&$\displaystyle{{}##}$\hfil
      \crcr#1\crcr}}\,}
\makeatother

\def\bull{\vrule height 1ex width 1ex depth -.1ex} 

\def\II{\mathaccent'27I} \def\1{\char'401}
\def\TT{\mathaccent'27T} \def\1{\char'401}
\def\DD{\mathaccent'27D} \def\1{\char'401}
\def\CC{\mathaccent'27C} \def\1{\char'401}
\def\BB{\mathaccent'27B} \def\1{\char'401}
\def\KK{\mathaccent'27K} \def\1{\char'401}
 \def\1{\char'401}
\def\YY{\mathaccent'27Y} \def\1{\char'401}

\title{The Renormalization Method and Quadratic-Like Maps}
\subjclass{Primary: 58F23, 54D05, Secondary: 30D05}
\thanks{Research is partially supported by an grant
from the NSF, by awards from the PSC-CUNY, and, at MSRI, by
NSF grant \# DMS-9022140.}
\author{Yunping Jiang}
\address{Department of Mathematics, Queens College of CUNY,
Flushing, NY 11367}
\email{yunqc@@yunping.math.qc.edu}

\begin{document}
\maketitle

\begin{abstract}
The renormalization of a quadratic-like map is studied. The
three-dimensional Yoccoz puzzle for an infinitely renormalizable
quadratic-like map is discussed. For an unbranched quadratic-like 
map having the {\sl a priori} complex bounds, the local connectivity of its 
Julia set is proved by using the three-dimensional Yoccoz puzzle. 
The generalized version of Sullivan's sector theorem is discussed 
and is used to prove his result that the Feigenbaum quadratic polynomial 
has the {\sl a priori} complex bounds and is unbranched. A dense subset on 
the boundary of the Mandelbrot set is constructed so that for every point 
of the subset, the corresponding quadratic polynomial is unbranched and 
has the {\sl a priori} complex bounds.
\end{abstract}

\vskip30pt
\centerline{\bf Contents}
\vskip10pt

\noindent 0. Introduction
\vskip5pt
\noindent 1. Quadratic Polynomials, Quadratic-Like Maps, and Hyperbolic
Geometry
\vskip5pt
\noindent 2. Renormalizable Quadratic-Like Maps
\vskip5pt
\noindent 3. Two-Dimensional Yoccoz Puzzles and Renormalizability
\vskip5pt
\noindent 4. Infinitely Renormalizable Quadratic Julia Sets and 
Three-Dimensional Yoccoz Puzzles
\vskip5pt
\noindent 5. A Generalized Sullivan's Sector Theorem
\vskip5pt
\noindent 6. Feigenbaum-Like Quadratic-Like Maps
\vskip5pt
\noindent 7. The Local Connectivity of Certain Infinitely Renormalizable 
Quadratic Julia Sets

\vfill

\eject
\section{Introduction}

Let $P_{c}(z) = z^{2}+c$ be a quadratic polynomial. A central problem 
in the study of the dynamics of $P_{c}$ is to understand the topology 
and geometry of the Julia set $J_{c}$ of $P_{c}$. The filled-in Julia 
set $K_{c}$ of $P_{c}$ is, by definition, the set of points not going 
to infinity under iterations of $P_{c}$. The Julia set $J_{c}$ of $P_{c}$
is, by definition, the boundary of $K_{c}$.

In order to have the more penetrating study of the dynamics of a quadratic
polynomial $P_{c}$, Douady and Hubbard~[DH3] introduced the concept of a
quadratic-like map. In \S 1, we review the definition of a
quadratic-like map and the work of Douady and Hubbard~[DH3] 
which proves that a quadratic-like map with connected  Julia set 
is hybrid equivalent to a unique quadratic polynomial. 
We also review the result due to Douady and Yoccoz (see~[HUB, MI2]) 
about the landing of external rays at a repelling periodic point of 
a quadratic polynomial. In the same section, we also review some fundamental 
results about the Julia set of a quadratic-like map and some basic facts 
of hyperbolic geometry.

The Mandelbrot set ${\cal M}$ is the set of complex parameters $c$ such that
the Julia set $J_{c}$ of $P_{c}$ is connected. The Julia set $J_{c}$ 
in the complement of ${\cal M}$ is a Cantor set. Douady and Hubbard~[DH1] 
proved that ${\cal M}$ is connected. They further conjectured that ${\cal
M}$ is locally connected. Many important research problems follow from this
conjecture. For example, the hyperbolicity conjecture (explained below)
would follow from this conjecture~[DH2]. A quadratic polynomial is hyperbolic if
it has an attractive or super-attractive periodic point in the complex
plane. Let ${\cal HP}$ be the set of parameters $c$ such that $P_{c}$ is
hyperbolic. The hyperbolicity conjecture says that ${\cal HP}$ is open and
dense in ${\cal M}$. 
To study the local connectivity of the Mandelbrot set ${\cal M}$, 
it would be helpful to answer the question: 
for which $c$ in ${\cal M}$ is the Julia set $J_{c}$ locally connected ?
Yoccoz made substantial progress in this direction. He proved that 
if $P_{c}$ is non-renormalizable (or finitely renormalizable), then the 
Julia set $J_{c}$ is locally connected. We discuss his result in \S 3. 
(Using this result, Yoccoz further proved that ${\cal M}$ is locally
connected at a finitely
renormalizable point $c$ (see~[HUB])). Also in \S 3, 
we discuss the two-dimensional Yoccoz puzzle of a quadratic-like 
map and its relation with the renormalizability of this quadratic-like map.

There remain many points in ${\cal M}$ which are infinitely renormalizable.
In \S 2, we define infinitely renormalizable quadratic-like maps and define 
infinitely renormalizable folding mappings, and discuss
the relation between the two definitions.

In \S 4, we prove one of our main results: for an unbranched 
infinitely renormalizable quadratic-like map having the {\sl a priori} complex
bounds, its Julia set is locally connected. We prove this result by using
the three-dimensional Yoccoz puzzle of an infinitely renormalizable
quadratic-like map.  Also in \S 4, we prove that 
the filled-in Julia set of any renormalization of a renormalizable
quadratic-like map, about the period of the two-dimensional Yoccoz puzzle, 
does not depend on the choices of renormalization domains. In particular, 
the renormalized filled-in Julia set is the limiting component in
two-dimensional Yoccoz puzzle containing the critical point.

In \S 5, we discuss Sullivan's sector theorem. We prove a generalized
version. The proof repeatedly applies
the sharpest version of Koebe's distortion theorem (see~[BIE]), and uses
Sullivan's idea about using hyperbolic contraction to trap points in a
hyperbolic neighborhood.

Using the generalized version of Sullivan's sector theorem, we 
prove in \S 6 his result that
the Feigenbaum quadratic polynomial has the {\sl a priori} complex
bounds and is unbranched. By combining this with the result in
\S 4, we complete the proof of the result which was first announced in~[JIH] 
and which says that the Julia set of the Feigenbaum quadratic polynomial 
is locally connected.

In \S 7, we construct a subset $\tilde{\Upsilon}$ of the Mandelbrot set
which is dense on the boundary $\partial {\cal M}$ of the Mandelbrot set 
${\cal M}$ such that for every point $c$ in this subset,
the corresponding quadratic polynomial $P_{c}(z)=z^{2}+c$ is unbranched and 
infinitely renormalizable and has the {\sl a priori} complex 
bounds. Thus, for $c$ in $\tilde{\Upsilon}$, the Julia set 
$J_{c}$ of a quadratic polynomial $P_{c}(z)=z^{2}+c$ is locally connected. 
A similar result concerning about the local 
connectivity of the Mandelbrot set at infinitely renormalizable points 
is proved in~[JI4].

This article benefited from reading J. Hubbard's
paper~[HUB] and J. Milnor's paper~[MI2], from talking with D. Sullivan, and
from listening the talk given by J.-C. Yoccoz in Denmark. During the writing of this article, 
several e-mails from C. McMullen provided great help. Some statements 
in this article are more precise because of his suggestions. 
Also from his e-mails, I decided to include the
papers~[JI1,JIH,JI2,JI3] into this self-contained article.
During my study of the work of
Sullivan, F. Gardiner and J. Hu provided a lot help. I thank them all.
This research was partially carried out~[JI1] when the author was 
in the IMS at Stony Brook. This research is partially supported by an grant
from the NSF and by awards from the PSC-CUNY. Research at
MSRI is supported in part by NSF grant \# DMS-9022140.

\section {Quadratic Polynomials, Quadratic-Like Maps, and
Hyperbolic Geometry}

Let $P_{c}(z)=z^{2}+c$ be a quadratic polynomial. Let $\bold C$ be
the complex plane. Let $V=\{ z \in \bold C\; |\;
|z| < r\}$ be a disk in $\bold C$. For $r$ large enough, $U=P^{-1}_{c}(V)$
is a simply connected domain, its closure is relatively compact in $V$,
and $P_{c}: U\rightarrow V$ is a holomorphic, proper map of degree two (see Fig.
1).
This is a model of an object defined by Douady and Hubbard~[DH3].

\vskip5pt
\centerline{\psfig{figure=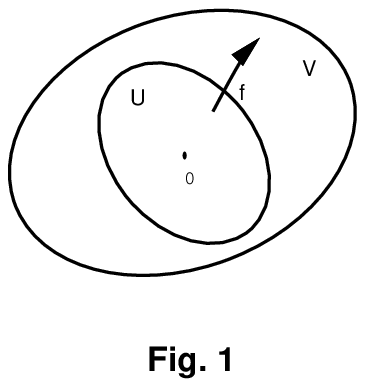}}
\vskip5pt

\proclaim Definition 1. A quadratic-like map is a triple
$(U,V,f)$ where $U$ and $V$ are simply connected domains isomorphic to a
disc with $\overline{U}\subset V$, and where $f:U\rightarrow V$ is
a holomorphic, proper map of degree two (see Fig. 1).

\vskip5pt
{\bf Remark 1.} A proper map $f$ means that $f^{-1}(K)$ is compact for
every compact set $K$. Douady and Hubbard~[DH3] also defined a
polynomial-like map.

\vskip5pt
Suppose $(U,V,f)$ is a quadratic-like map. We use
$$ K_{f}=\cap_{n=0}^{\infty} f^{-n}(U)$$
to denote the set of $z$ in $U$ such that images of $z$ under
iterations of $f$ are all in $U$. The set $K_{f}$ is called the filled-in
Julia set. It is a compact subset of $U$. The Julia set
$J_{f}$ of $f$ is the boundary of $K_{f}$.

A quadratic-like map $(U,V,f)$ has only one branched point $b$ at which
the
derivative $f'(b)$ of $f$ equals zero. We call $b$ the critical point of
$f$ and, without loss of generality, we always assume that $b=0$. The
following theorem is well-known.

\vskip5pt
\proclaim Theorem 1. Suppose $(U,V,f)$ is a quadratic-like map. If the
critical point
$0$ is not in $K_{f}$, then $K_{f}=J_{f}$ is a Cantor set of zero
Lebesgue measure in $\bold C$. And moreover, the set $K_{f}$ (or $J_{f}$) is
connected if and only if the critical point $0$ is in $K_{f}$.

\vskip5pt
Suppose that $\Omega$ is a domain in $\bold C$ and that $f$ is
a self-map of $\Omega$. A point $z$ in $\Omega$ is called a
periodic point of period $k\geq 1$ if $f^{\circ k}(z)=z$ and $f^{\circ i}(z) \neq
z$ for $1\leq i <k$. The number $\lambda =(f^{\circ k})'(z)$ is called the
multiplier (or eigenvalue) of $f$ at $z$. A periodic point of
period $1$ is called a fixed point.
A periodic point $z$ of $f$ is said to be repelling, attractive, or neutral if
$|\lambda| >1$, $0< |\lambda| <1$, or $|\lambda|=1$. Moreover, a
periodic point
is called super-attractive if $\lambda =0$, and is called parabolic if
$\lambda =e^{2\pi i p/q}$ where $p$ and $q$ are integers.
For an attractive or super-attractive or parabolic periodic point
of period $k\geq 1$ of a quadratic-like map $f: U\rightarrow V$, let
$O(p)=\{ f^{\circ i}(p)\}_{i=0}^{k-1}$ be the periodic orbit of $f$.
The set
$$
B(p) =\{ z\in V \; |\; f^{\circ n}(p) \rightarrow O(p) \hbox{ as }
n\rightarrow \infty\}
$$
is called the basin of $O(p)$. Let $CB(f^{\circ i}(p))$ be the connected
component containing $f^{\circ i}(p)$ of $B(p)$. The set $IB(p) =
\cup_{i=0}^{k-1} CB(f^{\circ i}(p))$ is called the immediate basin of
$O(p)$.
The proofs of the following two theorems can be found in Blanchard's
survey article~[BLA].

\vskip5pt
\proclaim Theorem 2. Suppose $z_{0}$ is a super-attractive periodic
point of period $k$
of a holomorphic function $f$ defined on $\Omega$. Then there is a
neighborhood $U$ of $z_{0}$, a unique holomorphic
diffeomorphism $h:U\rightarrow h(U)$ with $h(z_{0})=0$
and $h'(z_{0})=1$, and a unique integer $n>1$ such that
$$h\circ f^{\circ k}\circ h^{-1} (z) =z^{n}$$
on $h(U)$.

\vskip5pt
\proclaim Theorem 3. Let $(U,V,f)$ be a quadratic-like map and let
$J_{f}=\partial K_{f}$ be its Julia set. Let $E_{f}$ be the set of all
repelling periodic points of $f$. Then

\vskip3pt
\noindent {\bf (1)} $J_{f}$ is completely invariant, i.e., $f(J_{f})=J_{f}$ and
$f^{-1}(J_{f})=J_{f}$;

\vskip3pt
\noindent {\bf (2)} $J_{f}$ is perfect, i.e., $J_{f}'=J_{f}$, where $J_{f}'$
means the set of limit points of $J_{f}$;

\noindent {\bf (3)} $E_{f}$ is dense in the Julia set $J_{f}$, i.e., $\overline{E_{f}} =J_{f}$;

\vskip3pt
\noindent {\bf (4)} for any $z$ in $V$, the limit set of $\{ f^{-n}(z)
\}_{n=0}^{\infty}$ is $J_{f}$;

\vskip3pt
\noindent {\bf (5)} $J_{f}$ has no interior point;

\vskip3pt
\noindent {\bf (6)} If $f$ has an attractive or super-attractive or parabolic
periodic point $p$ in $V$, then the immediate
basin $IB(p)$ contains the critical point $0$ and
the critical orbit $O(0)=\{ f^{\circ n}(0)\}_{n=0}^{\infty}$;

\vskip3pt
\noindent {\bf (7)} If $f$ has neither any attractive, nor any super-attractive,
nor any neutral periodic points in $V$, then $K_{f}=J_{f}$.

\vskip5pt
Let $(U,V,f)$ and $(U',V',g)$ be two quadratic-like maps. They are
topologically conjugate if there is a homeomorphism $h$ from
a neighborhood $K_{f}\subset X \subset U$ to a neighborhood $K_{g}\subset
Y\subset U'$
such that $h\circ f=g\circ h$ on $X$ where $K_{f}$ and $K_{g}$ are filled-in
Julia sets. If $h$ is quasiconformal (see~[AH1]) (respectively,
holomorphic), then they are
quasiconformally (respectively, holomorphically) conjugate.
If $h$ can be chosen such that $h_{\overline{z}} =0$ a.e. on $K_{f}$,
then they are hybrid equivalent. Let
$$I(f)=\{ g \;\; |\; g \;\; \hbox{is hybrid equivalent to} \; f \}$$
be the inner class of $f$. The following theorem
is due to Douady and Hubbard.

\vskip5pt
\proclaim Theorem 4~[DH3]. If $(U,V,f)$ is a quadratic-like
map such that $K_{f}$ is connected, then there is a unique quadratic
polynomial $P(z) =z^{2} +c_{f}$ in $I(f)$.

\vskip5pt
Let $\overline{\bold C}=\bold C\cup \{ \infty\}$ be the extended complex
plane. Then $\infty$ is a super-attractive fixed point of any quadratic 
polynomial $P_{c}(z) =z^{2}+c$. The filled-in Julia set $K_{c}$ of $P_{c}$ 
is the set of all points not going to infinity under iterations of $P_{c}$.
Let ${\bf D}_{r}=\{ z\in \bold C\; |\; |z| <r\}$.
As we mentioned in the beginning of this section,
for $r>1$ large enough, $U=P^{-1}_{c}({\bf D}_{r})$ is a simply connected
domain and $\overline{U} \subset {\bf D}_{r}$.
Thus $(U,V,P_{c})$ for $V={\bf D}_{r}$ is a quadratic-like map. In
particular, $({\bf D}_{r},
{\bf D}_{r^{2}}, P_{0})$ is a quadratic-like map for every $r>1$; the filled-in
Julia set $K_{0}$ of $P_{0}$ is
the closed unit disk $\overline{\bf D}_{1}$.
Applying Theorem 2, there is a holomorphic diffeomorphism $h_{1}$ defined
on a neighborhood $\overline{\bold C}\setminus \overline{{\bf D}}_{r}$ (for
$r>1$ large) about $\infty$ such that $h_{1}(\infty)=\infty$,
$h_{1}'(\infty)=1$,
and such that $$ h_{1}^{-1}\circ P_{c}\circ h_{1} (z) =z^{2}$$
on $\overline{\bold C}\setminus \overline{{\bf D}}_{r}$. Let $B_{1}(\infty)
=h_{1}(\overline{\bold C}\setminus \overline{{\bf D}}_{r})$ and let
$B_{n}(\infty) =P^{-(n-1)}_{c}(B_{1}(\infty))$.
If the filled-in Julia set
$K_{c}$ of $P_{c}$ is connected, then all
$$P_{c}: B_{n}(\infty) \cap \bold C \rightarrow B_{n-1}(\infty)\cap \bold C$$
are unramified covering maps of degree two. We can inductively define
holomorphic diffeomorphisms $h_{n}$ on $\overline{\bold C}\setminus
\overline{{\bf D}}_{r^{1\over 2^{n}}}$ such that
$$ h_{n}^{-1}\circ P_{c}\circ h_{n} (z) =z^{2}$$
for $z$ in $\overline{\bold C}\setminus \overline{{\bf D}}_{r^{1\over
2^{n}}}$ and for $n>1$.
As $n$ tends to infinity, we get a holomorphic diffeomorphism
$h_{\infty}$ defined on $\overline{\bold C}\setminus \overline{{\bf D}}_{1}$
such that
$$ h_{\infty}^{-1}\circ P_{c}\circ h_{\infty} (z)   =z^{2} 
\eqno(1)$$
for all $z$ in $\overline{\bold C}\setminus \overline{{\bf D}}_{1}$. Therefore,
$B_{c}(\infty) =h_{\infty}(\overline{\bold C}\setminus \overline{\bf D}_{1})$ is
the basin of
$\infty$ for $P_{c}$ and $K_{c} =\overline{\bold C}\setminus B_{c}(\infty)$.
Furthermore, for every $r>1$ and for
$U_{r}=h_{\infty}({\bf D}_{r})$,
$(U_{r}, U_{r^{2}}, P_{c})$ is a quadratic-like map and its filled-in Julia set
is always $K_{c}$.
Let ${\bf S}^{R}=\{ z\in \bold C\; |\; |z|=R\}$
and let $s_{R}=h_{\infty}({\bf S}^{R})$ for $R>1$. Then
$$P_{c}(s_{R}) =s_{R^{2}}. \eqno(2)$$
The topological circle $s_{R}$ for every $R>1$ is
called an equipotential curve of $P_{c}$.
A curve
$$e_{\theta} =h_{\infty}(\{ z \in \overline{\bold C}\; |\; |z| >1,
\arg (z) =\theta\})$$
for $0\leq \theta <2\pi$ is called an external ray
of $P_{c}$. Then
$$P_{c}(e_{\theta})  = e_{2 \theta}.\eqno(3)$$

\vskip5pt
{\bf Remark 2.} Let
$$G (z) =\max \{ 0, \lim_{n\rightarrow \infty} {1\over 2^{n}} \log
|P_{c}^{\circ n}(z)|\}$$
be the Green's function of $K_{c}$ in $\overline{\bold C}$. Then $G(P_{c}(z))
=2G(z)$. For any $R>1$, the equipotential curve
$ s_{R} =G^{-1}(\log R)$ is a level curve of $G$.

\vskip5pt
If the mapping $h_{\infty}$ in Equation $(1)$ can be extended continuously to the unit circle
${\bf S}^{1}$, then we have a unique continuous map $H: \overline{\bf
C}\setminus {\bf D}_{1} \rightarrow \overline{\bold C}\setminus \KK_{f}$
such that $H|(\overline{\bf
C}\setminus \overline{\bf D}_{1}) =h_{\infty}$. Using $H$, we can define an
equivalence relation on ${\bf S}^{1}$: $z_{1}\sim z_{2}$
if and only if $H(z_{1})=H(z_{2})$. Let $[z]$ be
the equivalent class of $z$. Then $\tilde{P}_{0}([z]) =[P_{0}(z)]$
defines a map of the quotient space $X={\bf S}^{1}/\sim$, since
$z_{1}^{2} \sim z_{2}^{2}$ if $z_{1}\sim z_{2}$.
The dynamical system $(\tilde{P}_{0}, X)$ is
topologically conjugate to $(P_{c}, J_{c})$
by $\tilde{H}([z]) =H(z)$. The question arises:

\vskip5pt
\proclaim Question 1. For which $c$ can $h_{\infty}$ be extended
continuously to the unit circle ${\bf S}^{1}$ ?

\vskip5pt
A connected set $X$ in $\bold C$ is locally connected if for any point
$p$ in $X$ and any neighborhood $V$ about $p$, there is another
neighborhood $U\subset V$ about $p$ such that $U\cap X$ is connected.
The following classical theorem proved by Carath\'eodory
in one complex variable gives a sufficient and
necessary condition to extend $h_{\infty}$ continuously to ${\bf
S}^{1}$. The proof of this theorem can be found in~[MI1].

\vskip5pt
\proclaim Theorem 5~[CAR]. Let $h$ be a Riemann
mapping from ${\bf D}_{1}$ onto a simply
connected open domain $\Omega$. Then
$h$ can be extended continuously to the unit circle ${\bf S}^{1}$ if and
only if the boundary $\partial \Omega$ (as well as $\Omega$) is locally connected.

\vskip5pt
{\bf Remark 3.} If $\partial \Omega$ is a Jordan curve, then $h$ can be
extended to a homeomorphism from $\overline{\bf D}_{1}$ onto $\overline{\Omega}$.
Moreover, if $\partial \Omega$ is made of finite number of analytic curves, then the
extension restricted to the unit circle ${\bf S}^{1}$ has non-zero derivative
at every point other than a corner (see~[BIE]).

\vskip5pt
We have an equivalent question by just concerning the topology of a Julia
set now:

\vskip5pt
\proclaim Question 2. For which $c$ is $J_{c}$ locally connected ?

\vskip5pt
An external ray $e_{\theta}$ lands at $J_{c}$ if $e_{\theta}$ has only
one limit point at $J_{c}$. An external ray is periodic with period $m$
if $e_{\theta}\cap P_{c}^{\circ i}(e_{\theta}) =\{\infty \}$ for $1\leq
i <m$ and if $P_{c}^{\circ m}(e_{\theta}) =e_{\theta}$.
The following theorem is proved by Douady and Yoccoz (refer
to~[MI1,MI2,HUB]).

\vskip5pt
\proclaim Theorem 6. Let $P_{c}(z)=z^{2}+c$ be a quadratic
polynomial with connected Julia set $J_{c}$. Then every
repelling periodic point of $P_{c}$ is a landing point of finitely
many periodic external rays with the same period.

\vskip5pt
Let ${\cal S}$ be a Riemann surface and let $(\tilde{{\cal S}}, \pi)$
be the universal cover of ${\cal S}$, where $\tilde{{\cal S}}$ is a simple
connected Riemann surface and $\pi: \tilde{{\cal S}}
\rightarrow {\cal S}$ is the universal covering map.
From the Uniformization Theorem (see~[AH2]), we can identify
$\tilde{{\cal S}}$
with one of the extended complex plane $\overline{\bold C}$, the complex
plane $\bold C$, or the open unit disk ${\bf D}_{1}$.
A Riemann surface ${\cal S}$ is hyperbolic if $\tilde{{\cal S}}={\bf D}_{1}$.

Let ${\cal D}$ be the hyperbolic disk which is the open unit disk with the
hyperbolic metric
$$d_{H}s = {|dz|\over 1-|z|^{2}}.$$
Let $d_{H}$ be the hyperbolic distance.
Every holomorphic diffeomorphism $h: {\bf D}_{1}\rightarrow {\bf D}_{1}$,
which is not
a linear fractional transformation, strictly decreases the hyperbolic
distance $d_{H}$, i.e,
$$ d_{H} (h(z_{1}), h(z_{2})) < d_{H}(z_{1}, z_{2})$$
for all $z_{1}$ and $z_{2}$ in ${\cal D}$. For a hyperbolic Riemann surface
${\cal S}$, one can define the hyperbolic distance $d_{H,S}$ from $d_{H}$
and $\pi$.
Any holomorphic map $h: {\cal S} \rightarrow {\cal S}$,
which is not an isometry with respect to $d_{H,{\cal S}}$,
strictly decreases this hyperbolic distance $d_{H, S}$,
i.e., for any $z_{1}$ and $z_{2}$ in ${\cal S}$,
$$d_{H,{\cal S}}(h(z_{1}), h(z_{2})) < d_{H,{\cal S}}(z_{1}, z_{2}).$$

A bounded domain $\Omega$ in $\bold C$ is a hyperbolic Riemann
surface.
An important family of hyperbolic Riemann surfaces is the family of
bounded doubly connected
domains in $\bold C$. A bounded doubly connected domain $\Omega$ is
called an
annulus. From complex analysis (see~[BIE]), any
annulus is holomorphically diffeomorphic to a unique round annulus
$A_{r}=\{ z\in \bold C \; |\; r < |z| <1\}$ for $0< r<1$.
The number
$$\hbox{mod}(\Omega) =- \log r$$
is called the modulus of $\Omega$. It is a
conformal invariant, i.e., $\hbox{mod}(h(\Omega)) =\hbox{mod}(\Omega)$ whenever
$h$ is
a conformal homeomorphism from
$\Omega$ onto $h(\Omega)$.

Let $\Omega$ be a bounded doubly connected domain in $\bold C$. Then the
complement $\bold C\setminus \Omega$ of $\Omega$ in $\bold C$ has two components.
One is a connected and simply connected compact set $E$ and the other
is an unbounded
set $F$. From the Gr\"otzsch argument (see~[AH1]),
the bounded component $E$ is a single point if and
only if $\hbox{mod}(\Omega)=\infty$ (see~[BRH]).

Let $\Omega$ be an annulus. If $\Omega_{1} \subseteq \Omega$ is
a subannulus, then
$$\hbox{mod}(\Omega_{1}) \leq \hbox{mod}(\Omega).
\eqno(4)$$
If $\Omega_{1}$, $\Omega_{2} \subseteq \Omega$ are two disjoint
subannuli, then
$$ \hbox{mod}(\Omega_{1}) +\hbox{mod}(\Omega_{2}) \leq
\hbox{mod}(\Omega). \eqno(5)$$
The proofs of these two inequalities can be found in the book of Ahlfors~[AH1].

Let $E$ be a connected and simply connected compact subset of the open
unit disk ${\bf D}_{1}$. Let
$$ \hbox{mod}({\bf D}_{1}, E) = \sup_{\Omega} \{ \hbox{mod}(\Omega),
\hbox{ where } \Omega \subseteq {\bf D}_{1}\setminus E \hbox{ is a round subannulus}\, \}.$$
Note that $\hbox{mod}({\bf D}_{1}, E) =\infty$ if $E$ is a single point.

Let $\hbox{diam}_{H}(E)
=\sup_{z_{1}, z_{2} \in E} d_{H}(z_{1}, z_{2})$ be the hyperbolic diameter
of $E$ in ${\bf D}_{1}$.
The following theorem can be found in the book of McMullen~[MC1].

\vskip5pt
\proclaim Theorem 7. The hyperbolic diameter $d_{H}(E)$ and the
$\hbox{mod}({\bf D}_{1}, E)$ are inversely related:
$$ \hbox{diam}_{H}(E) \rightarrow 0 \qquad \iff \qquad \hbox{mod}({\bf D}_{1},
E) \rightarrow \infty$$
and
$$ \hbox{diam}_{H}(E) \rightarrow \infty \qquad \iff \qquad
\hbox{mod}({\bf D}_{1}, E) \rightarrow 0.$$
More precisely, there is a constant $C>0$ such that
$$ C^{-1} \hbox{diam}_{H} (E) \leq \exp (-\hbox{mod}({\bf D}_{1}, E)) \leq
C\hbox{diam}_{H}(E)$$
when $\hbox{diam}_{H}(E)$ is small, while
$$ {C \over \hbox{diam}_{H} (E)} \geq \hbox{mod}({\bf D}_{1}, E) \geq
C^{-1}\exp (-\hbox{diam}_{H}(E))$$
when $\hbox{dima}_{H}(E)$ is large.

\section {Renormalizable Quadratic-Like Maps}

Let $(U_{0},V_{0},f_{0})$ be a quadratic-like map and suppose
its filled-in Julia set $K_{f_{0}}$ is connected.
The map $f_{0}: U_{0}\rightarrow V_{0}$ is renormalizable if there are an integer $n\geq 2$ and
an open subdomain $U_{1}$ of $U_{0}$ such that $0\in U_{1}$ and such that
$(U_{1}, V_{1}, f_{0}^{\circ n})$
is a quadratic-like map
with connected filled-in Julia set, where $V_{1}=f^{\circ n}(U_{1})$.
Let $f_{1}=f_{0}^{\circ n}|U_{1}$.
The filled-in Julia set $K_{f_{1}} = K_{f_{1}}(n, U_{1})$
is {\sl a priorily} dependent on the choice of $U_{1}$ (refer to Theorem 13).
The domain $U_{1}$ is called a renormalization and
$(U_{0},V_{0},f_{0})$ is
called renormalizable about $n$. Otherwise, $(U_{0},V_{0},f_{0})$ is called
non-renormalizable.

The quadratic-like map $(U_{0},V_{0},f_{0})$ is infinitely renormalizable
if there is a strictly increasing sequence of integers $\{ n_{k}\}_{k=1}^{\infty}$
such that $(U_{0},V_{0},f_{0})$ is renormalizable about $n_{k}$ for
$k=1$, $2$, $\ldots$. Otherwise, $(U_{0},V_{0},f_{0})$ is called finitely
renormalizable.

A quadratic-like map $(U_{0},V_{0},f_{0})$ is called real if $f_{0}(U_{0}\cap \bold R^{1})
\subset V_{0}\cap \bold R^{1}$ and $g=f_{0}|U_{0}\cap \bold R^{1}$ is a
real folding map. Let $(U_{0}, V_{0}, f_{0})$ be a real
quadratic-like map and suppose its filled-in Julia set $K_{f_{0}}$ is connected.
Suppose $g=f_{0}|U_{0}\cap \bold R^{1}$ has a fixed point $p \in \bold R^{1}$ with positive
multiplier $f_{0}'(p)$. Let $p'\neq p$ be another inverse image of $p$ under
$g$, that is, $g(p')=p$. Conjugating by a linear fractional transformation, we
may assume that $p=-1$ and $p'=1$ and that
$$[-1,1]=\cap_{n=0}^{\infty}g^{-n}(V_{0}\cap \bold R^{1})=K_{f}\cap {\bf
R}^{1}.$$
Hence $g$ is a folding map with unique quadratic critical point $0$.

Let $g$ be a folding map of $[-1,1]$ such that $g(-1)=g(1)=-1$ and such that
$0$ is a unique quadratic critical point. We say $g$ is renormalizable
about $n>1$ if there is a subinterval $I$ of $[-1,1]$ such that
$0\in \II$, such that $g^{\circ
i}(I)\cap \II=\emptyset$ for all $0< i<n$, and such that $g^{\circ n}(I)\subseteq
I$. Otherwise, $g$ is non-renormalizable. We say $g$ is infinitely
renormalizable if there is a strictly increasing sequence $\{
n_{k}\}_{k=1}^{\infty}$ such that $g$ is renormalizable about $n_{k}$ for
all $k>0$. Otherwise, $g$ is called finitely renormalizable.
The next theorem
shows that for a real quadratic-like map,
both definitions of renormalization are essentially equivalent.

\vskip5pt
\proclaim Theorem 8. Let $(U_{0},V_{0},f_{0})$ be a real
quadratic-like
map and suppose its filled-in Julia set $K_{f_{0}}$ is connected. Suppose $f_{0}$ has neither
neutral, nor attractive, nor super-attractive periodic points.
Then $f_{0}$ is renormalizable if and only if the folding map $g=f_{0}|[-1,1]$ is renormalizable.

\vskip5pt
{\bf Proof.} Suppose the folding map $g=f_{0}|[-1, 1]$ is renormalizable. This
means that there is a maximal closed subinterval $I$ of $[-1,1]$ and an integer $n>1$
such that $0$ is in $\II$, such that $g^{\circ i}(I)\cap \II =\emptyset$ for $0<i< n$, and
such that $g^{\circ n}(I)\subseteq I$. One of the endpoints of $I$ is fixed by
$g^{\circ
n}$. It is a repelling fixed point. Take a neighborhood $T'$ of $I$ such that
$f^{\circ n}|(L\cup R)$ is expanding, where $L\cup R=T'\setminus I$.
Let $I_{i}=g^{\circ i}(I)$ for $0\leq i\leq n$.
The inverse $h_{n}$ of $g^{\circ (n-1)}:I_{1}\rightarrow I_{n}$ is a diffeomorphism and can be
extended to a diffeomorphism on an open interval $T \supset I$.
Take $T\subset T'$.
Let $M'=h_{n}(T)$ and $M=g^{-1}(M')$. Then $\overline{M}$ is a subset of $T$.
Because the critical orbit $CO =\{ f_{0}^{\circ k}(0)\}_{k=0}^{\infty}$ of
$f_{0}$ is in the real line $\bold R^{1}$, $h_{n}$ can be extended to
$V_{1}= (V_{0}\setminus \bold R^{1})\cup
T$ analytically. The image $V_{1}'$ of $V_{1}$ under this extension is
contained in $U'=(U_{0}\setminus \bold R^{1})\cup M'$. Let $U_{1}
=f_{0}^{-1}(V_{1}')$. Then $\overline{U}_{1} \subset V_{1}$ and $(U_{1},
V_{1}, f_{1})$ for $f_{1} =f_{0}^{\circ n}|U_{1}$ is a quadratic-like map.
Since $g^{\circ n}(I)\subseteq I$, the filled-in Julia set $K_{f_{1}}$ is
connected. Therefore, $(U_{0}, V_{0}, f_{0})$ is renormalizable about $n$ and
$U_{1}$ is a renormalization.

Suppose $(U_{0}, V_{0}, f_{0})$ is renormalizable about $n>1$.
Let $U_{1}$ be a renormalization and set $V_{1}=f^{\circ
n}(U_{1})$.
Then $(U_{1}, V_{1}, f_{1})$ for $f_{1}=f_{0}^{\circ n}|U_{1}$ is a
quadratic-like map with the connected filled-in
Julia set $K_{f_{1}}$. Let $I=K_{f_{1}} \cap \bold R^{1}$.
Then $I=\cap_{i=0}^{\infty} f^{- i}_{1}(V_{1}\cap \bold R^{1})$.
For every $1\leq i<n$, $f^{\circ i}(I)\cap \II =\emptyset$ else
$f^{\circ n}$ would have at least three fixed points in $I\cup f^{\circ
i}(I)$ with one
of them either attractive or parabolic. Since $f_{1}(0)$ is in $I$,
$f_{0}^{\circ n}(I)\subseteq I$. Therefore, $g=f_{0}|[-1, 1]$ is a renormalizable
folding map.
\hfill \bull

\section {Two-Dimensional Yoccoz Puzzles and Renormalizability}

In this section, we discuss a technique in the study of
non-renormalizable quadratic polynomials, due to Yoccoz, and some of its
applications to renormalization.

Let $P_{c}(z)=z^{2}+c$ be a quadratic polynomial with connected
filled-in Julia set $K_{c}$.
The external ray $e_{0}$ of $P_{c}$ is the only one fixed by
$P_{c}$ (see Fig. 2). It lands either at a repelling or
at a parabolic fixed point
$\beta$ of $P_{c}$ (see~[MI1]). Suppose $\beta$ is
repelling. Applying Theorem 6, we see that $e_{0}$ is the
only external ray landing at $\beta$. Thus $K_{c}\setminus \{ \beta \}$
is connected. We call $\beta$ the non-separate fixed point of
$P_{c}$.
Let $\alpha\neq \beta$ be the other fixed point of $P_{c}$. If
$\alpha$ is either an attractive
or a super-attractive fixed point, then $J_{c}=K_{c}\setminus
(\cup_{n=0}^{\infty} P^{-n}_{c} \big( D(\alpha)\big)$ for a small
disk centered at
$\alpha$. The Julia set $J_{c}$ is a Jordon curve; every external ray
lands
at a unique point in $J_{c}$ (see Remark 3). If $\alpha$
is a repelling fixed point,
there are at least two periodic external rays landing at $\alpha$.
We use $R_{0}(\alpha)$ to denote the union of a cycle of periodic
external rays of period $q$ landing at
$\alpha$ (see Fig. 2). The set $R_{0}(\alpha)$ cuts $\bold C$ into
finitely many simply
connected domains $\Omega_{0}$, $\Omega_{1}$, $\ldots$, $\Omega_{q-1}$.
Each domain contains points in the Julia set $J_{c}$.
Thus $K_{c}\setminus \{ \alpha\}$ is disconnected. We call $\alpha$ the
separate fixed point of $P_{c}$.

\vskip10pt
\centerline{\psfig{figure=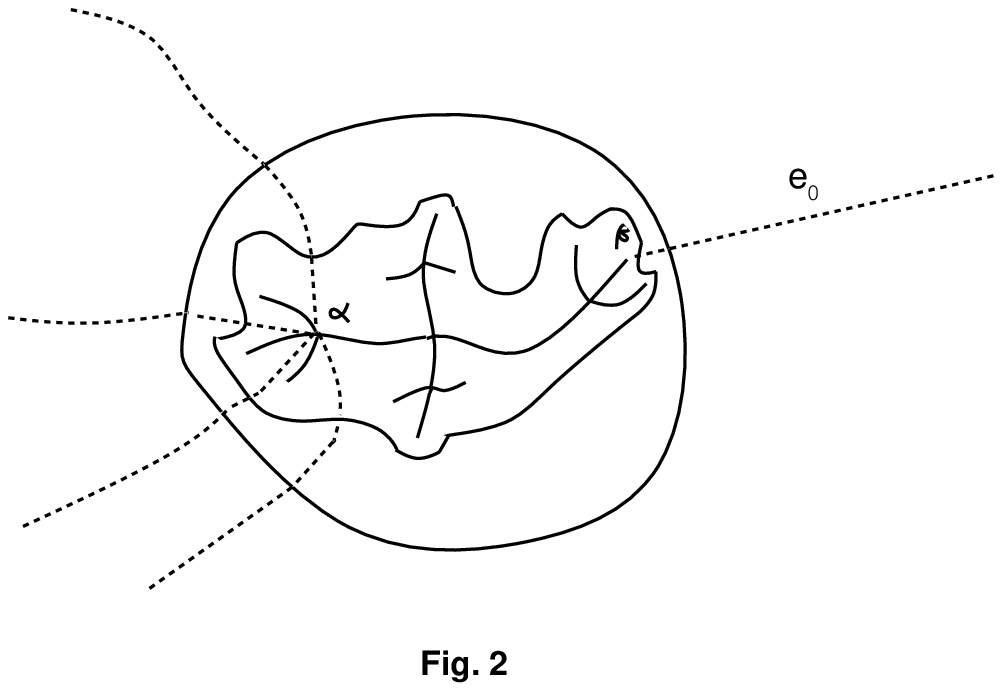}}
\vskip5pt

Henceforth, we assume that the fixed points $\beta$ and $\alpha$
are both repelling. Let
$s_{r}$ be a fixed equipotential curve of $P_{c}$ and let $U_{r}$ be the
open domain bounded by $s_{r}$. Then $(U_{\sqrt{r}}, U_{r}, P_{c})$ is a
quadratic-like map. The set $R_{0}(\alpha)$ cuts
$U_{r}$ into finitely many simply connected domains.
Let $C_{0}$ be the closure of the domain containing $0$, and let $B_{0,i}$ be the closure of the domain
containing $P_{c}^{\circ i}(0)$ for $1\leq i <q$.
Since $R_{0}(\alpha)$ is forward invariant under
$P_{c}$, the image under $P_{c}$ of $C_{0}\cap K_{c}$ or
$B_{0,i}\cap K_{c}$, for every $1\leq i <q$,
is the union of some of $C_{0}\cap K_{c}$, $B_{0,1}\cap K_{c}$, $\ldots$,
$B_{0,q-1}\cap K_{c}$.
The set $\eta_{0}=\{ C_{0}, B_{0,1}, \ldots, B_{0,q-1}\}$
is called the original partition. We note that it is not a
Markov partition because $P_{c}|C_{0}$ is a proper, holomorphic map
of degree two.
(But $P_{c}|B_{0,i}$ is a holomorphic diffeomorphism for every $1\leq i <q$.)

Let $\Gamma_{n}=P_{c}^{-n}(\alpha)$ and let
$R_{n}(\alpha) =P_{c}^{-n}\big( R_{0}(\alpha)\big)$ for $n\geq 0$.
The set $R_{n}(\alpha)$ is the union of some external rays
landing at points in $\Gamma_{n}$; it cuts the domain $U_{r^{1\over 2^{n}}}$
into a finite number of simply connected domains.
Let $C_{n}$ be the closure of the domain containing $0$ and let $B_{n, 1}$,
$\ldots$, $B_{n,
k_{n}}$ be the closures of others. Since $P_{c}(\Gamma_{n}) =\Gamma_{n-1}$,
the image of $C_{n}$
or $B_{n,i}$ under $P_{c}$ is one of
$C_{n-1}$, $B_{n-1,1}$, $\ldots$,
$B_{n-1,k_{n-1}}$, for $1\leq i\leq k_{n}$ and $n\geq 1$.
Then $P_{c}|C_{n}$ is
holomorphic, proper branch covering map of degree two; all
$P_{c}|B_{n,i}$ are holomorphic diffeomorphisms.
The set $\eta_{n}=\{ C_{n}, B_{n,1}, \ldots, B_{n,k_{n}}\}$
is called the $n^{th}$-partition. The sequence $\eta=\{ \eta_{n}
\}_{n=0}^{\infty}$ is called the two-dimensional Yoccoz puzzle for
$P_{c}$. A similar puzzle
for certain cubic polynomial is constructed by Branner and Hubbard~[BRH].
Yoccoz used this puzzle
while studying the local connectivity of a non-renormalizable
quadratic polynomial as follows.

Let $\Gamma_{\infty} =\cup_{n=0}^{\infty}\Gamma_{n}$.
For any $x$ in $K_{c}\setminus \Gamma_{\infty}$, there is one and only one
sequence $\{ D_{n}(x)\}_{n=0}^{\infty}$ such that $x\in D_{n}(x)\in \eta_{n}$.
For any $x$ in $\Gamma_{\infty}$, there are $q$ such sequences.
We call such a sequence
$$ x\in \cdots \subseteq D_{n}(x)\subseteq D_{n-1}(x) \subseteq \cdots
\subseteq D_{1}(x)\subseteq D_{0}(x)$$
an $x$-end. In particular,
$$ 0\in \cdots \subseteq C_{n}\subseteq C_{n-1} \subseteq \cdots
\subseteq C_{2}\subseteq C_{1} \subseteq C_{0}$$
is called the critical end.

\vskip5pt
Suppose $D_{n+1}\subseteq D_{n}$ are domains in $\eta_{n+1}$ and
$\eta_{n}$ for $n\geq 0$. Define $A_{n}=D_{n}\setminus \DD_{n+1}$.
If $D_{n+1} \subset \DD_{n}$, then $\AA_{n}$ is a non-degenerate
annulus and
its modulus (denoted $\hbox{mod}(A_{n})$) is greater than $0$.
Otherwise, $A_{n}$
is a degenerate annulus and its modulus $\hbox{mod}(A_{n})$ is zero.
The domain $D_{n}$ is critical if $D_{n}=C_{n}$ and
$D_{n+1}=C_{n+1}$; it is
semi-critical if $D_{n}= C_{n}$ but $D_{n+1}\neq C_{n+1}$; it is
non-critical if $D_{n}\neq C_{n}$.
If $D_{n}$ is non-critical, then
$$\hbox{mod}(A_{n})=\hbox{mod}(P_{c}(A_{n}))$$
since $P_{c}: A_{n}\rightarrow P_{c}(A_{n})$ is a conformal homeomorphism.
If $D_{n}$ is critical, then
$$\hbox{mod}(A_{n})={\hbox{mod}(P_{c}\big( A_{n})\big)\over 2}$$
since $P_{c}:A_{n}\rightarrow P_{c}(A_{n})$ is a
proper, holomorphic, unramified covering map of degree two. If $D_{n}$ is semi-critical,
then $${\hbox{mod}(P_{c}\big( A_{n})\big)\over 2} \leq \hbox{mod}(A_{n}).$$

For a point $x$ in $K_{c}$ and an $x$-end
$$ x\in \cdots \subseteq D_{n}(x)\subseteq D_{n-1}(x) \subseteq \cdots
\subseteq D_{1}(x)\subseteq D_{0}(x),$$
we define
$$D_{nm}(x) =P^{\circ m}_{c}(D_{n+m}(x)) (= D_{n}(P^{\circ m}_{c}(x)))$$
and
$$A_{nm}(x) =D_{nm}(x) \setminus \DD_{(n+1)m}(x)$$
for $n,m\geq 0$.
The tableau $T(x)=(a_{nm})_{n\geq 0, m\geq 0}$ is an $\infty \times
\infty$-matrix defined as follows:
$a_{nm}=1$ if $D_{nm}(x)=C_{n}$, and $a_{nm}=0$
if $D_{nm}(x) \neq C_{n}$.
The tableau $T(0)=\Big( a_{nm}^{0} \Big)_{n\geq 0,m\geq 0}$
is called the critical tableau. Note that
$$ P_{c}^{\circ m}(x) \in \cdots \subseteq D_{nm}(x) \subseteq
D_{(n-1)m}(x) \subseteq \cdots \subseteq D_{1m}(x) \subseteq D_{0m}(x)$$
is a $P_{c}^{\circ m}(x)$-end.

\vskip5pt
\proclaim Lemma 1.
The tableau $T(x)$ satisfies the following rules:
\vskip5pt
{\bf (T1)} if $a_{nm}=1$ for $n,m\geq 0$, then $a_{im}=1$ for all $0\leq i\leq n$,
\vskip5pt
{\bf (T2)} if $a_{nm}=1$ for $n,m\geq 0$, then $a_{(n-i)(m+j)} =a_{(n-i)j}^{0}$ for all
$0\leq i+j\leq n$,
\vskip5pt
{\bf (T3)} if
\vskip3pt
$\hskip10pt$ {\bf i)} $a_{nm}=1$ and $a_{(n+1)m}=0$ for $n,m\geq 0$ and if
\vskip3pt
$\hskip10pt$ {\bf ii)} $a_{(n-i)(m+i)}=1$ for $1\leq i\leq n$ and $a_{(n-j)(m+j)} =0$
for $0< j<i$,
\vskip3pt
\noindent then $a_{(n-i+1)i}^{0}=1$ implies $a_{(n-i+1)(m+i)} =0$.

\vskip5pt
{\bf Proof.} Rule {\bf (T1)} is valid because
$$ 0\in \cdots \subseteq C_{n}\subseteq C_{n-1}\subseteq \cdots
\subseteq C_{1}\subseteq C_{0}.$$
Rule {\bf (T2)} is valid because if
$a_{nm}=1$, then $D_{im}(x)=C_{i}$ for all $0\leq i\leq n$. Thus
$D_{(n-i)(m+j)}(x) = f^{\circ j}(C_{n-i})$ for all $0\leq i+j\leq n$.

We prove Rule {\bf (T3)}. Conditions {\bf i)} and {\bf ii)} imply
that $P_{c}^{\circ i}: C_{n} \rightarrow C_{n-i}$ is a degree two
holomorphic proper branch covering map. Condition {\bf i)}
and Rule {\bf (T2)} imply $a_{k(m+j)}=a_{kj}^{0}$ for $0\leq k+j \leq n$.
Now $a_{(n-i+1)i}^{0}=1$ and $a_{j0}^{0}=1$ for all $0\leq j<\infty$ imply
that $P_{c}^{\circ i}: C_{n+1} \rightarrow C_{n-i+1}$ is a degree two
holomorphic proper
branch covering map. Thus $C_{n+1} =P_{c}^{-i}(C_{n-i+1})\cap C_{n}$.
Assume $a_{(n-i+1)(m+i)} =1$. Then $D_{(n+1)m}(x) =
P_{c}^{-i}(C_{n-i+1})\cap C_{n} =C_{n+1}$. This contradicts
$a_{(n+1)m}=0$.
\hfill \bull

\vskip5pt
\proclaim Lemma 2. For any domain $D$ in $\eta_{n}$, for $n\geq 0$,
$D\cap K_{c}$ is connected.

\vskip5pt
{\bf Proof.} Since the domain $D$ is bounded by finitely many
external rays $\Pi=\{ e_{\theta_{i}} \}_{i=1}^{m}$ and by some
equipotential curve, then $\partial D \cap K_{c}$
consists of a finite number of points $\{ p_{i}\}_{i=1}^{n}$. Every
$p_{i}$ is a landing point of two
external rays in $\Pi$. Suppose $D\cap K_{c}$ is not connected for
some $D$ in $\eta_{n}$. Then there are two disjoint open sets $U$ and
$V$ such that $D\cap K_{c} =(D\cap K_{c} \cap U)
\cup (D\cap K_{c}\cap V)$. Suppose that $p_{1}$, $\ldots$, $p_{k}$ are in $U$
and that $p_{k+1}$, $\ldots$, $p_{n}$ are in $V$. The two external rays in
$\Pi$ landing at $p_{i}$ cut $\bold C$ into two open domains.
Let $W_{i}$ be the one which is disjoint from $D$.
Then $U'=U\cup \cup_{i=1}^{k}W_{i}$ and $V'=V\cup \cup_{i=k+1}^{n}W_{i}$
are two disjoint open sets and $K_{c} =(U'\cap K_{c}) \cup (V'\cap
K_{c})$. This contradicts the fact that $K_{c}$ is connected.
\hfill \bull

\vskip5pt
Lemma 2 implies that for any $x$-end
$$ x\in \cdots \subseteq D_{n}(x)\subseteq D_{n-1}(x) \subseteq \cdots
\subseteq D_{1}(x)\subseteq D_{0}(x),$$
the intersection $L_{x} =\cap_{n=0}^{\infty} D_{n}(x)$
is a compact connected non-empty set containing $x$. Let
$T(x)=(a_{nm})_{n\geq 0, m\geq 0}$ be the tableau of the $x$-end. It
is non-recurrent if there is an
integer $N\geq 0$ such that $a_{nm}=0$ for all $n\geq N$ and all $m\geq
1$. Otherwise, $T(x)$ is recurrent.

\vskip5pt
\proclaim Lemma 3. If $T(x)$ is non-recurrent, then $L_{x} =\{ x\}$.

\vskip5pt
{\bf Proof.} Suppose $N\geq 0$ is an integer such that
$a_{nm}=0$ for all $n\geq N$ and
all $m\geq 1$. Then, for $n >N$, every
$$P_{c}^{\circ (n-N-1)}: D_{n1}(x) \rightarrow
D_{N(n-N)}(x)$$
is a holomorphic diffeomorphism. Thus $\hbox{mod}(A_{n0}(x))$ is
greater than or equal to
$\hbox{mod}(A_{N(n-N)}(x))/2$
for every $n>N$. There are only finitely many different annuli in $\{
A_{Nm}(x)\}_{m=0}^{\infty}$ because $\eta_{N}$ has only finitely many
domains.
If there were infinitely many non-degenerate annuli in
$\{A_{Nm}(x)\}_{m=0}^{\infty}$, then there would be infinitely many
non-degenerate annuli in $\{ A_{n0}(x)\}_{n=0}^{\infty}$ whose moduli are
the same.
This would imply that
$$\hbox{mod} (D_{0}(x)\setminus L_{x}) \geq \sum_{n=0}^{\infty}
\hbox{mod}(A_{n0}(x)) =\infty .$$
Therefore $L_{x} =\{ x\}$.

If there are only finitely many non-degenerate annuli in
$\{A_{Nm}(x)\}_{m=0}^{\infty}$, The proof uses results from hyperbolic
geometry. Let $B_{N1}$ be the domain in $\eta_{N}$
containing
the critical value $P_{c}(0)$. Since $P_{c}^{\circ i}(x)$
do not enter $C_{N+1}$ for all $0< i<
\infty$, then $P_{c}^{\circ i}(x)$ does not enter $B_{N1}$ for all $2\leq
i<\infty$.
Let us thicken $B_{Ni}$ to an open simply
connected domain $\tilde{B}_{Ni}$ such that $B_{Ni} \subset \tilde{B}_{Ni}$
and such that $P_{c}(0)$ is not in $\tilde{B}_{Ni}$ for $1< i\leq k_{N}$.
The map $P_{c}$ has two inverse
branches $g_{i1}$ and $g_{i2}$ defined on $\tilde{B}_{Ni}$ for every
$1< i\leq k_{N}$. we consider $\tilde{B}_{Ni}$ to be a hyperbolic
Riemann surface with hyperbolic distance $d_{H,i}$ for every $1< i\leq
k_{N}$, where $k_{N}$ is the number of elements in $\eta_{N}$. Then if $g_{ik}$,
for $1< i\leq k_{N}$ and
$k=1$ or $2$, sends $B_{Ni}$ into $B_{Nj}$ for some $1< j\leq k_{N}$,
then it strictly contracts
these hyperbolic distances; more precisely, there is a constant
$0< \lambda <1$ such
that $d_{H,j}(g_{ik}(x), g_{ik}(y)) < \lambda d_{H,i}(x,y)$ for
$x$ and
$y$ in $B_{Ni}$ and for $k=0$ and $1$. Therefore, there is a constant $C>0$ such that
for any $D_{n1}(x)$ and for any $n>N$,
$$d(D_{n1}(x)) =\max_{y,z\in D_{n0}(x)} |y-z| \leq C \lambda^{n-N-1}$$
since $D_{(n-i)i}(x)$ is in one of $B_{N2}$, $\ldots$,
$B_{Nk_{N}}$ for
every $2\leq i\leq n-N$.
Thus $d(D_{n0}(x))$ tends to zero as $n$ goes to infinity and
$L_{x} =\{x\}$.
\hfill \bull

\vskip5pt
The critical end
$$ 0\in \cdots \subseteq C_{n}\subseteq C_{n-1} \subseteq \cdots
\subseteq C_{2}\subseteq C_{1} \subseteq C_{0}$$
is important. Let $A_{n0}(0)=C_{n}\setminus \CC_{n+1}$.

\vskip5pt
\proclaim Lemma 4. If $\sum_{n=0}^{\infty} \hbox{mod}(A_{n0}(0)) =\infty$, then for any $x$ in
$K_{c}$ and any $x$-end
$$ x\in \cdots \subseteq D_{n}(x)\subseteq D_{n-1}(x) \subseteq \cdots
\subseteq D_{1}(x)\subseteq D_{0}(x),$$
$L_{x} =\{ x\}$.

\vskip5pt
{\bf Proof.} Consider the tableau $T(x)=(a_{nm})_{n\geq 0, m\geq 0}$. If
$T(x)$ is non-recurrent, the lemma follows from Lemma 3.

Suppose $T(x)$ is recurrent. If there is a column which is entirely
$1$'s, then there are integers $M\geq 0$ and $N\geq 0$
such that $a_{iM}=1$ for all $i\geq 0$ and $a_{nm}=0$ for all
$n\geq N$ and $0\leq m< M$. Thus
$P_{c}^{\circ M}: D_{n0}(x) \rightarrow D_{(n-M)M}=C_{n-M}$ is a holomorphic
diffeomorphism for every $n\geq N$. This implies
$$ m(D_{0}(x)\setminus L_{x}) \geq \sum_{n=0}^{\infty} \hbox{mod} (A_{n0}(x))
\geq \sum_{n=N}^{\infty} \hbox{mod} (A_{n0}(x)) =\sum_{n=N-M}^{\infty}
\hbox{mod} (A_{n0}(0)) =\infty.$$
So $L_{x}=\{ x\}$.

Suppose that there is no column which is entirely $1$'s.
Let $N>0$ be an integer such that $a_{n0}=0$ for $n\geq N$.
For any $n\geq N$, let $m_{n}>0$ be the integer such that $a_{nm_{n}}=1$
and $a_{ni}=0$ for $0\leq i<m_{n}$. Then
$P_{c}^{\circ m_{n}}:
D_{(n+m_{n}-1)0}(x) \rightarrow D_{(n-1)m_{n}}(x)=C_{n-1}$ is a
holomorphic diffeomorphism. Remember that $A_{(n-1)0}(x) =D_{(n-1)0}(x)\setminus
\DD_{n0}(x)$
and $A_{(n-1)0}(0) =C_{n-1}\setminus \CC_{n}$. We have
$$ \hbox{mod}\big( A_{(n+m_{n}-1)0}(x)\big) =\hbox{mod}\big(
A_{(n-1)0}(0)\big) .$$
Let $q_{n}=n+m_{n}-1$. Then $q_{N}<q_{N+1}<\cdots < q_{n}
<q_{n+1}<\cdots$. Thus
$$ \eqalign{\hbox{mod} (D_{0}(x)\setminus L_{x}) &\geq \sum_{n=N}^{\infty}
\hbox{mod}\big( A_{n0}(x)\big) \geq \sum_{n=N}^{\infty}
\hbox{mod}\big( A_{q_{n}0}(x)\big) \cr&=\sum_{n=N}^{\infty} \hbox{mod}\big( A_{(n-1)0}(0)\big) =\infty.}$$
This implies that $L_{x}=\{x\}$.
\hfill \bull

\vskip5pt
The first column of the critical tableau $T(0)=(a_{nm}^{0})_{n\geq 0, m
\geq 0}$ is entirely $1$'s.
If $T(0)$ has another column which is entirely $1$'s, that is, if there
is an integer $m>0$ such that $a_{im}=1$ for all $i\geq 0$, then
we call $T(0)$ a periodic critical tableau.

\vskip5pt
\proclaim Theorem 9~[YOC].
The critical tableau $T(0)$ is periodic if and only
if $P_{c}$ is renormalizable.

\vskip5pt
{\bf Proof.} Suppose $T(0)$ is periodic. Let $n_{1}>0$ be the smallest integer such that
$a_{in_{1}}=1$ for all $i\geq 0$. Let $N\geq 0$ be the smallest integer
such that $a_{ij}=0$ for all $i\geq N$ and $0<j<n_{1}$. For any $n\geq
n_{1}+N$, $P_{c}^{\circ n_{1}}:C_{n} \rightarrow C_{n-n_{1}}$ is a degree
two proper holomorphic branch covering map.
Thus $\{ P^{\circ
kn_{1}}(0)\}_{k=0}^{\infty}$ is contained in $C_{n_{1}+N}$.
If $C_{n_{1}+N}\subset \CC_{N}$, then
$P_{c}^{\circ n_{1}}:\CC_{n_{1}+N}\rightarrow \CC_{N}$ is a
quadratic-like map with connected filled-in Julia set and
$\CC_{n_{1}+N}$ is a renormalization about $n_{1}$. Thus $P_{c}$ is
renormalizable.

In general, let us consider a small open disk $D(\alpha)$ centered at
the separate fixed point $\alpha$ of $P_{c}$ such that
$$\overline{D}(\alpha) \subset D'(\alpha)=P_{c}^{\circ
n_{1}}\big( D(\alpha)\big)$$ and such that
$$D'(\alpha) \cap \{ P_{c}^{\circ i}(0)\}_{i=0}^{n_{1}+N}=\emptyset .$$
Thicken $C_{0}$ and $B_{0i}$ for $1\leq i\leq q-1$ as follows.
Suppose $C_{0}$ (respectively, $B_{0i}$) is bounded by two
external rays $R_{\theta_{1}}$ and $R_{\theta_{2}}$ of angles
$\theta_{1}$
and $\theta_{2}$. Let $\epsilon >0$ be a small number such that the
domains
$$U_{1}=\cup_{\theta_{1}-\epsilon <\theta <\theta_{1}+\epsilon}
(R_{\theta}\setminus D_{\alpha})\qquad \hbox{and}\qquad
U_{2}=\cup_{\theta_{2}-\epsilon <\theta <\theta_{2}+\epsilon}
(R_{\theta}\setminus D_{\alpha})$$
are disjoint from $\{ P_{c}^{\circ
i}(0)\}_{i=0}^{n_{1}+N}$. Let
$$\tilde{C}_{0} = (U_{1}\cup C_{0}\cup U_{2}\cup D(\alpha))\cap U_{r}
\qquad (\hbox{respectively,}
\qquad \tilde{B}_{0i} = (U_{1}\cup B_{0i}\cup U_{2}\cup D(\alpha))\cap
U_{r} )$$
where $U_{r}$ is the domain bounded by the equipotential curve $s_{r}$.
Let
$$\tilde{\eta}_{0} =\{ \tilde{C}_{0}, \tilde{B}_{01}, \ldots ,
\tilde{B}_{0(q-1)}\}$$
and let
$$\tilde{\eta}_{n} =P^{-n}(\tilde{\eta}_{n}) = \{ \tilde{C}_{n}, \tilde{B}_{n1},
\ldots , \tilde{B}_{nk_{n}}\}$$
for $1\leq n \leq n_{1}+N$.
The diffeomorphism $g=P_{c}^{-(n_{1}-1)}:
C_{N}\rightarrow P_{c}(C_{n_{1}+N})$ can be extended to
$\tilde{C}_{N}$. Let $B'$ be the image of $\tilde{C}_{N}$ under $g$. Then
$\tilde{C}_{n_{1}+N}=P_{c}^{-1}(B')$.
Let $\tilde{\CC}_{n}$ denote the interior of $\tilde{C}_{n}$ for
$0\leq n\leq n_{1}+N$. Then
$\tilde{C}_{n_{1}+N}\subset \tilde{\CC}_{N}$.
Thus
$$P_{c}^{\circ n_{1}}: \tilde{\CC}_{n_{1}+N} \rightarrow \tilde{\CC}_{N}$$
is a quadratic-like map and
$\tilde{\CC}_{n_{1}+N}$ is a renormalization about $n_{1}$.
This proves the ``only if'' part.

Now suppose $P_{c}$ is renormalizable. Let $U_{1}$ be a renormalization
about $n_{1}$,
that is, $(U_{1}, V_{1}, f_{1})$ is a quadratic-like map with connected
filled-in Julia set $K_{f_{1}}$ where $f_{1}=P_{c}^{\circ n_{1}}|U_{1}$ and
$V_{1}=f_{1}(U_{1})$. The map $f_{1}$ has two fixed points $\beta_{1}$ and
$\alpha_{1}$ in $U_{1}$. Let
$$ \alpha_{1}\in \cdots \subseteq D_{n}(\alpha_{1}) \subseteq
D_{n-1}(\alpha_{1}) \subseteq \cdots \subseteq D_{1}(\alpha_{1}) \subseteq
D_{0}(\alpha_{1})$$
be an $\alpha_{1}$-end. There is a $D_{k}(\alpha_{1})$ such that $K_{f_{1}}
\subset D_{k}(\alpha_{1})$ and $D_{k}(\alpha_{1}) = C_{k}$.
Since
$$K_{f_{1}} \subseteq f_{1}^{-1} (U_{1}\cap C_{k}) \subseteq U_{1}\cap C_{k},$$
then $P^{\circ n_{1}}_{c}$ sends $C_{k+in_{1}}$ to
$C_{k+(i-1)n_{1}}$ for all $i>0$. Thus $T(0)$ is periodic. It is the ``if'' part.
\hfill \bull

\vskip5pt
We define a function $\tau$ on the set $\bold N$ of natural numbers by
using the critical
tableau $T(0) =(a_{nm}^{0})_{n\geq 0, m\geq 0}$ as follows: $\tau (n)=m$
if $a_{(n-i)i}^{0}=0$ for $0 <i<n-m$ and if $a_{m(n-m)}^{0}=1$; if there
is no such integer $m>0$, then $\tau(n)=-1$.

If the critical tableau $T(0)$ is periodic, then there are integers $n_{1}>0$ and $N\geq 0$ such
that $a_{in_{1}}=1$ for all $i\geq 0$ and such that $a_{ij}=0$ for $i\geq N$ and
$0< j< n_{1}$. Thus $\tau (n) =n-n_{1}$ for $n\geq N+n_{1}$.

If the critical tableau $T(0)$ is non-recurrent, then there is the smallest
integer $N \geq 0$ such that $a_{nm}=0$ for all $n\geq N$ and $m>0$.
Thus the image $\tau (\bold N)$ is contained in the finite set
$\{ -1, 0, 1, \ldots, N-1\}$.

If the critical tableau $T(0)$ is not periodic and is recurrent, then
every row of $T(0)$ has infinitely many $1$'s and every column (except
for the $0^{th}$-column) has a $0$. An integer $n\geq 0$
is noble if for every entry $a_{nk}^{0}$ such that $a_{nk}^{0}=1$,
we have $a_{(n+1)k}^{0}=1$.

\vskip5pt
\proclaim Lemma 5. If the critical tableau $T(0)$ is not
periodic and is recurrent,
then the function $\tau$ satisfies the following properties:

\vskip3pt
{\bf (i)} For any integer $m\geq 0$, $\tau^{-1}(m)$ is not empty.

\vskip5pt
{\bf (ii)} If $m\geq 0$ is noble, then $\tau^{-1}(m)$ contains at least
two different integers.

\vskip5pt
{\bf (iii)} If $\tau (n)=m$ and if $m$ is noble, then $n$ is also noble.

\vskip5pt
{\bf (iv)} If $\tau^{-1}(m)$ contains only one integer $n$, then $n$ is
noble.

\vskip5pt
{\bf Proof.} We prove {\bf (i)} first. Consider any $m$-row in
$T(0)$ for $m\geq
0$. Let $k>0$ be the integer such that $a_{mi}^{0}=0$ for $0<i<k$ and such that
$a_{mk}^{0}=1$. From {\bf (T1)}, $a_{(m+k-i)i}^{0} =0$ for $0< i <k$.
Thus $\tau (m+k)= m$.

To prove {\bf (ii)}, suppose $m\geq 0$ is noble. Let $k$ be the
same integer
as that in the proof of {\bf (i)}. Let $m_{1}$ be the integer such that
$a_{m_{1}k}^{0}=1$ and such that $a_{ik}^{0}=0$ for all $i>m_{1}$. Consider
$a_{(m_{1}-k)(2k)}^{0}$, $a_{(m_{1}-2k)(3k)}^{0}$, $\ldots$, and
$a_{(m_{1}-(i-1)k)(ik)}^{0}$ where
$m_{1}-ik\leq m< m_{1}-(i-1)k$. From the tableau
rules {\bf (T1)} and {\bf (T3)},
$a_{(m_{1}-k+1)(2k)}^{0}=0$, $a_{(m_{1}-2k+1)(3k)}^{0}=0$, $\ldots$,
$a_{(m_{1}-(i-1)k+1)(ik)}^{0}=0$.
If $m=m_{1}-ik$, from the tableau rules {\bf (T1)} and {\bf (T3)},
$a_{m(ik)}^{0}=1$. Since $m$ is noble, $a_{(m+1)(ik)}^{0}=1$. But from the tableau
rules {\bf (T1)} and {\bf (T3)}, $a_{(m+1)(ik)}^{0} =0$. The contradiction
implies that $m>m_{1}-ik$.

Now from the tableau rule {\bf (T2)}, $a_{m(k+m_{1}-m)}^{0}=0$.
Let $k_{1}>k+m_{1}-m$ be the integer such that
$a_{mi}^{0}=0$ for $k+m_{1}-m < i < k_{1}$ and such that $a_{mk_{1}}^{0}=1$.
Then $a_{(m+k_{1}-i)i}^{0} =0$ for $0<i< k_{1}$ and $a_{mk_{1}}^{0}=1$. This says
that $\tau (m+k_{1})=m$.

To prove {\bf (iii)}, suppose $\tau(n)=m$ where $m$ is noble.
For any $a_{nk}^{0}=1$, since $a_{(n-i)i}^{0}=0$ for $0<i<n-m$ and
$a_{m(n-m)}^{0}=1$ and since the tableau rule {\bf (T1)}, we have
$a_{(n-i)(k+i)}^{0}=0$ for $0<i< n-m$ and $a_{m(k+n-m)}^{0}=1$. Since
$m$ is noble, $a_{(m+1)(n-m)}^{0}=1$. Assume $a_{(n+1)k}^{0}=0$. From
the tableau rule {\bf (T3)}, $a_{(m+1)(k+n-m)}^{0}=0$. This
contradicts to that $m$ is noble. Thus $a_{(n+1)k}^{0}=1$. This means that
$n$ is noble.

Now we prove {\bf (iv)}.
Suppose $n>0$ is the only integer such that $\tau (n) =m$.
We first consider $a_{(m+1)(n-m)}^{0}$. If $a_{(m+1)(n-m)}^{0}=0$, then we
would have an integer $k> n-m$ such that $a_{mi}^{0}=0$ for $n-m< i<k$
and such that $a_{mk}^{0}=1$. From the tableau rule {\bf (T1)},
$a_{(m+k-i)i}^{0}=0$ for $0 < i < k$. This would imply
that $\tau (m+k)=m$, which contradicts the assumption.
Thus, $a_{(m+1)(n-m)}^{0}=1$.
If there is an entry $a_{nk_{1}}^{0}=1$ with
$a_{(n+1)k_{1}}^{0}=0$, then $k_{1}> n-m$ and, from the tableau
rule
{\bf (T3)}, $a_{m(n-m+k_{1})}^{0}=1$ and $a_{(m+1)(n-m+k_{1})}^{0}=0$.
Consider the smallest integer $k_{2}>n-m+k_{1}$ such that $a_{mi}^{0}=0$ for
$k_{1}+n-m < i<k_{2}$
and $a_{mk_{2}}^{0}=1$. From the tableau rule {\bf (T1)},
$a_{(m+k_{2}-i)i}^{0}=0$ for $k_{1}+n-m < i <
k_{2}$. So we can find another integer $n_{0} \geq k_{2}-k_{1} +m >n$
such that $\tau (n_{0}) =m$. This would contradict the assumption.
\hfill \bull

\vskip5pt
\proclaim Theorem 10~[YOC].
Suppose $P_{c}(z)=z^{2}+c$ is a
recurrent quadratic polynomial.
The critical tableau $T(0)$ is periodic
if and only if $L_{0}$ contains more than one point.

\vskip5pt
{\bf Proof.}
We use the same notation as in the proof of Theorem 9 and the proof of
Lemma 5.
Suppose $T(0)$ is periodic.
Then for $n> N+n_{1}$
$$P_{c}^{\circ n_{1}} : C_{n+1}\rightarrow C_{n-n_{1}+1}$$
is a degree two proper holomorphic branch covering map.
Replacing $C_{n+1}$ by
$\tilde{C}_{n+1}$ if it is necessary, we may assume that this map is a
quadratic-like map. Since $L_{0}$ is the filled-in Julia set of this
map, it contains more than one point.
This is the ``only if'' part.

To prove the ``if'' part, suppose $T(0)$ is not periodic.
We will prove that $L_{0}$ contains only
one point. Since $P_{c}$ is recurrent, there are infinitely many $1$'s
in every row of $T(0)$, that is, $T(0)$ is recurrent.
Consider the first partition
$$\eta_{1}=\{ C_{1}, B_{11}, \ldots, B_{1(q-1)}, B_{01},\ldots B_{0(q-1)}\},$$
where $B_{0i}=B_{0,i}$ for $1\leq i <q$ and where $B_{1i}\subseteq C_{0}$ and
$P_{c}(B_{1i})=B_{0,i}$ for $1\leq i <q$. (Remember that $\eta_{0}=\{ C_{0},
B_{0,1}, \ldots , B_{o, q-1}\}$ is the original partition.)
Let $c(n)= P_{c}^{\circ
n}(0)$. If the critical orbit
$CO=\{ c(n) \}_{n=0}^{\infty}$
is contained in the union
$$C_{1}\cup B_{01}\cup
\ldots \cup B_{0(q-1)},$$
then $T(0)$ is periodic of period $q$. Hence
there must be one critical value $c(n)$ in $B_{11}\cup \cdots
\cup B_{1(q-1)}$.
Let $c(n)$ be in $B_{1i}$. The annulus $A_{0n}(0) = C_{0}\setminus
\BB_{1i}$
is non-degenerate.  Pull back $A_{0n}(0)$ by $P_{c}$ along
$A_{i(n-i)}(0)$ for $0\leq i\leq n$; we get a
non-degenerate annulus $A_{n0}(0)$.

Now consider $\tau^{-k}(n)$. For each $m$ in $\tau^{-k}(n)$,
$$\hbox{mod}\big( A_{m0}(0)\big) \geq {\hbox{mod}\big( A_{n0}(0)\big) \over 2^{k}}.$$
If the number of $\tau^{-k}(n)$ is greater than or equal to $2^{k}$ for
every $k>0$,
then
$$ \eqalign{\hbox{mod}(C_{0}\setminus L_{0}) \geq \sum_{m=1}^{\infty}
\hbox{mod}\big( A_{m0}(0)\big) &\geq \sum_{k=1}^{\infty}
\sum_{m\in \tau^{-k}(n)} \hbox{mod} \big( A_{m0}(0)\big)
\geq \sum_{k=1}^{\infty}
\hbox{mod} \big( A_{n0}(0)\big) \cr
&=\Big( \hbox{mod}\big( A_{n0}(0)\big)
\Big)\cdot \sum_{k=1}^{\infty}
 1  = \infty.\cr}$$
So $L_{0} =\{0\}$.

If there is an integer $k>0$ such that the number of $\tau^{-k}(n)$ is
less
than $2^{k}$, then there are pre-images $m>q$ of $n$ under iterates of
$\tau$ such that $m$ is the only pre-image of $q$ under $\tau$. From
{\bf (iii)}, $m$ is noble. Hence $\tau^{-k}(m)$ are all noble and
contain at least $2^{k}$ different integers. Moreover
$$\hbox{mod}(A_{p0}(0)) = {\hbox{mod}(A_{m0}(0))\over 2^{k}}$$
for every $p$ in $\tau^{-k}(m)$.
Therefore,
$$ \eqalign{\hbox{mod}(C_{0}\setminus L_{0}) \geq
\sum_{k=1}^{\infty} \hbox{mod}\big( A_{k0}(0)\big) &\geq
\sum_{k=1}^{\infty}
\sum_{p\in \tau^{-k}(m)} \hbox{mod} \big( A_{p0}(0)\big)
\geq \sum_{k=1}^{\infty}
\hbox{mod} \big( A_{m0}(0)\big) \cr
&= \Big( \hbox{mod}\big(
A_{m0}(0)\big) \Big)\cdot \sum_{k=1}^{\infty}
 1  = \infty.\cr}$$
Again we have $L_{0}=\{ 0\}$. This completes the ``if'' part.
\hfill \bull

\vskip5pt
\proclaim Theorem 11~[YOC]. If $P_{c} (z) =z^{2}+c$ is a
non-recurrent or recurrent
non-renorm\-alizable quadratic polynomial, then its filled-in Julia set $K_{c}$
is locally connected.

\vskip5pt
{\bf Proof.} Let $\alpha$ be the separate fixed point of $P_{c}$.
Construct
the two-dimensional Yoccoz puzzle for $P_{c}$. For any $x$ in $K_{c}$,
let
$$ x\in \cdots \subseteq D_{n}(x) \subseteq D_{n-1}(x) \subseteq \cdots
\subseteq D_{1}(x) \subseteq D_{0}(x)$$
be an $x$-end. If $P_{c}$ is non-recurrent, then the critical
tableau is non-recurrent. Lemma 3 and Lemma 4 imply that
the diameter
$d(D_{n}(x))$ tends to zero as $n$ goes to infinity.
If $P_{c}$ is recurrent and non-renormalizable, then $T(0)$ is recurrent
and is not periodic.
Lemma 4 and Theorem 10 imply that the diameter
$d(D_{n}(x))$ tends to zero as $n$ goes to infinity.

If $x$ is not a preimage of $\alpha$ under any iterate of $P_{c}$, then
$x$ is an
interior point of $D_{n}(x)$ for all $n\geq 0$.
From Lemma 2, $\{D_{n}(x)\}$ is a basis of connected neighborhoods at $x$.
If $x$ is a preimage of $\alpha$ under some iterate of $P_{c}$, then
there are $q$ different $x$-ends,
$$ x\in \cdots \subseteq D_{i,n}(x) \subseteq D_{i,(n-1)}(x) \subseteq
\cdots \subseteq D_{i,1}(x) \subseteq D_{i,0}(x)$$
where $q$ is the period of the external rays landing at $\alpha$.
Let $\tilde{D}_{n}(x) =\cup_{i=1}^{q} D_{i,n}(x)$. Then $x$ is an interior
point of $\tilde{D}_{n}$. Since $K_{c}\cap D_{1,n}(x)$, $\ldots$, $K_{c}\cap
D_{q,n}(x)$ have a common point $x$, from Lemma 2, $K_{c}\cap
\tilde{D}_{n} (x)$ is connected.
So $\{ \tilde{D}_{n}(x) \}_{n=0}^{\infty}$ is a basis of connected
neighborhoods at $x$.
\hfill \bull

\vskip5pt
From Theorem 4, all arguments in this section apply
to a quadratic-like
map. Suppose that $(U,V,f)$ is a quadratic-like map and that its filled-in
Julia set $K_{f}$ is connected. Suppose two fixed points
$\beta$ and
$\alpha$ of $f$ are repelling.
Let $\beta$ be the non-separate fixed point of $f$, that is,
$K_{f}\setminus \{\beta\}$ is connected, and let $\alpha$ be the
separate
fixed point of $f$, that is, $K_{f}\setminus \{ \alpha\}$ is
disconnected.
Since $(U,V,f)$ is hybrid equivalent to a quadratic polynomial
$P_{c}$, there is a quasiconformal homeomorphism $H$ defined on $V$ such
that
$$H\circ f=P_{c}\circ H$$
on $U$. We call $e_{\theta, f}=H^{-1}(e_{\theta}\cap H(U))$ the external
ray of angle $\theta$ of $f$ where $e_{\theta}$ is the external ray of
$P_{c}$ of angle $\theta$.

Two points $H(\beta)$ and $H(\alpha)$ are non-separate and separate
fixed points of $P_{c}$, respectively. Suppose $\Gamma$ is the union of
a cycle of periodic
external rays landing at $H(\alpha)$. Let $\Gamma'=H^{-1}(\Gamma \cap
H(U))$. The set $\Gamma'$ cuts the domain $U$ into $q$ domains. Each of
them contains points in the filled-in Julia set $K_{c}$. Let $C_{0}$ be the
domain containing $0$ and let $B_{0, i}$ be the domain containing $f^{\circ
i}(0)$ for $1\leq i< q$. The partition
$$\eta_{0}=\{ C_{0}, B_{0, 1}, \ldots , B_{0, q-1}\}$$
is called the original partition for $f$. Let $\Gamma_{n}'=f^{-n}(\Gamma')$
and $U_{n}=f^{-n}(U)$.
Then $\Gamma_{n}'$ cuts $U_{n}$ into finitely many domains. Let $C_{n}$ be
the domain containing $0$ and $B_{n, i}$ for $1\leq i\leq k_{n}$ be others.
Then
$$\eta_{n}=\{ C_{n}, B_{n, 1}, \ldots , B_{n, k_{n}}\}$$
is called the $n^{th}$-partition for $f$. We use $f^{-n}(\eta_{0})$ to
denote $\eta_{n}$, i.e., $\eta_{n}=f^{-n}(\eta_{0})$, for $1\leq n< \infty$.
We have that
$f(C_{n})$ and $f(B_{n,i})$ for $1\leq i\leq k_{n}$ are in $\eta_{n-1}$ for
$n>0$ (set $k_{0}=q-1$).  We call $\eta=\{ \eta_{n}\}_{n=0}^{\infty}$
the two-dimensional Yoccoz puzzle of $f$.
Let $\Lambda =\cup_{n=0}^{\infty} f^{-n}(\alpha)$. Let
$L_{0}=\cap_{n=0}^{\infty} C_{n}$ be the
connected component of $K_{f}\setminus \Lambda$ containing $0$. We state
Theorems 10 and 11 in the following form.

\vskip5pt
\proclaim Theorem 12~[YOC]. Suppose $(U,V,f)$ is a
recurrent quadratic-like map. Then $(U, V, f)$ is
renormalizable if and only if $L_{0}$ contains more than one point.
Moreover, if $(U,V,f)$ is non-renormalizable, then any connected component of
$K_{f}\setminus \Lambda$ consists of only one point and $K_{f}$ is locally
connected.
\hfill \bull

\section {Infinitely Renormalizable Quadratic Julia Sets and
Three-Dimensional Yoccoz Puzzles}

Suppose $(U,V,f)$ is a renormalizable quadratic-like map with connected filled-in
Julia set $K_{f}$. Let $\eta=\{ \eta_{n}\}_{n=0}^{\infty}$ be
the two-dimensional Yoccoz puzzle for $f$. From the previous
section, the critical tableau $T(0)=(a_{nm}^{0})_{n\geq 0, m\geq 0}$ is periodic of period $n_{1}$.
Let
$$0\in \cdots \subseteq C_{n}\subseteq C_{n-1}\subseteq \cdots \subseteq
C_{1}\subseteq C_{0}$$
be the critical end. There is an integer $N>0$ such that $a_{ij}^{0} = 0$ for all $i\geq N$ and
for all $0< j< n_{1}$.
Let $f_{0}=f^{\circ n_{1}}|\CC_{N+n_{1}}$. Then
$f_{0}: \CC_{N+n_{1}}
\rightarrow \CC_{N}$ is a proper, holomorphic branch cover of degree two.
Assume $C_{N+n_{1}} \subset \CC_{N}$. (Otherwise, we can replace $C_{n}$
with $\tilde{C}_{n}$ (see Theorem 9)).
Then $(\CC_{N+n_{1}}, \CC_{N}, f_{0})$ is a quadratic-like map and its
filled-in Julia set is $L_{0} =\cap_{n=0}^{\infty} C_{n}$.

\vskip5pt
\proclaim Theorem 13. Suppose $(U,V,f)$ is a renormalizable
quadratic-like map with connected filled-in Julia set $K_{f}$. For any
renormalization $U_{1}$ about $n_{1}$, let $f_{1}=f^{\circ n_{1}}|U_{1}$
and let $V_{1}=f_{1}(U_{1})$. Then the filled-in Julia set $K_{f_{1}}$ (or
the Julia set $J_{f_{1}}$) of $(U_{1}, V_{1}, f_{1})$ is always $L_{0}$ (or
$\partial L_{0}$).

\vskip5pt
{\bf Proof.} Let $U'=\CC_{N+n_{1}}\cap U_{1}$ and let $U''$ be the connected
component of $U'$ containing $0$. Let $f_{2}=f^{\circ n_{1}}|U''$ and
$V''=f_{2}(U'')$. Then $f_{2} : U'' \rightarrow V''$
is a degree two branch covering. It is also proper because
$$f_{2}^{-1}(K)=f^{-1}_{0}(K)\cap f^{-1}_{1}(K)$$
for any compact set $K$ of
$V''$. Since $C_{N+n_{1}} \subset \CC_{N}$ and $\overline{U}_{1} \subset
V_{1}$, then
$$C_{N+n_{1}} \cap \overline{U}_{1} \subset \CC_{N}\cap V_{1}.$$
Thus $\overline{U''} \subset V''$ because they are the connected
components of $U'$
and $\CC_{N}\cap V_{1}$ containing $0$. Both $U''$ and $V''$ are simply connected and
isomorphic to a disc because they are intersections of
simply connected domains each of which is isomorphic to a disc.
Therefore, $(U'', V'', f_{2})$ is a quadratic-like
map. Let $K_{f_{2}}$ be the filled-in Julia set of $(U'', V'', f_{2})$.
Since a filled-in Julia set is completely invariant and since $0$ is in $U''$,
the two inverse images of $0$ under $f_{2}$ are in
$K_{f_{2}}$. But these two points are also inverse images of $0$ under
$f_{0}$ and under $f_{1}$. Therefore, they are both in $L_{0}$ and
in $K_{f_{1}}$. Using this argument, the set $\Xi$ of all inverse images of
$0$ under iterates of $f_{2}$ is contained in $K_{f_{2}}$ and is also
contained in $L_{0}$
and in $K_{f_{1}}$. Therefore,
$$K_{f_{1}} =K_{f_{2}}=L_{0} \hbox{ (or } J_{f_{1}}=J_{f_{2}}=\partial L_{0} \hbox{)}$$
because each of $\partial L_{0}$, $J_{f_{1}}=\partial K_{f_{1}}$, and
$J_{f_{2}}=\partial K_{f_{2}}$ is the closure of $\Xi$ (see Theorem 3).
\hfill \bull

\vskip5pt
As we saw in \S 2, the definition of the filled-in Julia set of
a renormalization about $n_{1}$ depends on choices of domains in
renormalization.
But the renormalized filled-in Julia set is actually canonical; it is
independent of choices of domains in
renormalization and is the limiting piece containing $0$ in the
two-dimensional
Yoccoz puzzle from Theorem 12.

Suppose $(U_{1},V_{1},f_{1})$ is a recurrent renormalizable
quadratic-like map with connected filled-in Julia set $K_{1}$.
We call $K_{1}$ (or $J_{1}$) a quadratic
filled-in Julia set (or quadratic Julia set). It is renormalizable if the corresponding
quadratic-like map is renormalizable.
Let $\beta_{1}$ and $\alpha_{1}$ be the
non-separate and separate fixed points of $f_{1}$, i.e.,
$K_{1}\setminus \{ \beta_{1}\}$ is still connected and
$K_{1}\setminus
\{\alpha_{1}\}$ is disconnected. Let $\Lambda_{1}=\cup_{n=0}^{\infty}
f_{1}^{-n}(\alpha_{1})$. Let $K_{2}=L_{0}$ be the connected component of
$K_{1}\setminus \Lambda_{1}$ containing $0$.
From Theorem 12, $K_{1}$ is renormalizable if and only if
$K_{2}$ contains more than one point. The quadratic filled-in Julia set
$K_{2}$ is called the renormalization of $K_{1}$.

Inductively, let $K_{i}$ be the renormalization of $K_{i-1}$. Let
$f_{i}=f_{i-1}^{\circ n_{i-1}}$ for $i\geq 2$, where $n_{i-1}$ is the period of the
critical tableau $T^{i-1}(0)= (a_{nm}^{0}(i-1))_{n\geq 0, m \geq 0}$ of the
two-dimensional Yoccoz puzzle for $(U_{i-1}, V_{i-1}, f_{i-1})$.  Let
$\beta_{i}$ and $\alpha_{i}$ be the
non-separate and the separate fixed points of $f_{i}$, i.e.,
$K_{i}\setminus \{ \beta_{i}\}$ is still connected and $K_{i}\setminus
\{\alpha_{i}\}$ is disconnected. Let $\Lambda_{i}=\cup_{n=0}^{\infty}
f_{i}^{-n}(\alpha_{i})$ and let $K_{i+1}$ be the connected component of
$K_{i}\setminus \Lambda_{i}$ containing $0$. Then
$K_{i}$ is renormalizable if and only if
$K_{i+1}$ contains more than one point. Here $K_{i}$, for $i>1$, is called
the $i^{th}$-renormalization of $K_{1}$.
Theorem 12 can be
generalized as follows.

\vskip5pt
\proclaim Theorem 14~[YOC]. Suppose that $(U_{1}, V_{1}, f_{1})$ is a
recurrent quadratic-like map and that $K_{1}$ is its filled-in
Julia set. The quadratic Julia set
$K_{1}$ is finitely renormalizable if and only if 
there is an integer $m\geq 1$ such that
$K_{1}$, $\ldots$, $K_{m}$ contains more than one point 
and such that $K_{m+1}$ contains only the point $0$. 
Moreover, if $K_{1}$ is finitely renormalizable, 
then $K_{1}$ is locally connected.

\vskip5pt
{\bf Proof.} The first part of the theorem follows directly from Theorem 12.
We prove the second part.
Let $\alpha_{m}$ be the separate fixed point of $f_{m}$.
Let $\Gamma_{m}$ be a cycle of periodic external rays of $f_{1}$
landing at $\alpha_{m}$ (refer to the end of \S 3).
Using $\Gamma_{m}$, we can construct the two-dimensional Yoccoz puzzle:
let $\eta^{m}_{0}$ be the set consisting of the closures of
the connected components of $V_{1}\setminus \Gamma_{m}$. Let
$\eta_{n}^{m}=f_{1}^{-n}(\eta_{0}^{m})$ for $n\geq 1$. Let
$C^{m}_{n}$ be the member of $\eta^{m}_{n}$
containing $0$. Since $f_{m}$ is non-renormalizable, we use a proof similar
to that of Theorem 10 to show that $\sum_{n=0}^{\infty}
\hbox{mod}(A_{n0}^{m}(0))
=\infty$, where $A_{n0}^{m}(0)=C_{n}^{m}\setminus \CC_{n+1}^{m}$.
Applying Lemma 4,
for every $x$-end
$$ x\in \cdots \subseteq D^{m}_{n}(x) \subseteq D^{m}_{n-1}(x) \subseteq
\cdots \subseteq D^{m}_{1}(x)
\subseteq D^{m}_{0}(x),$$
$L_{x}=\cap_{n=0}^{\infty} D^{m}_{n}(x)$ contains only $x$.
By using a similar argument to the proof of Theorem 11, we
can now show that $K_{1}$ is locally connected.
\hfill \bull

\vskip5pt
Now let us consider an infinitely renormalizable
quadratic-like map $(U_{1}, V_{1}, f_{1})$.
Let $K_{1}$ be the filled-in Julia set of $f_{1}$. Let $K_{i}$ be the
$i^{th}$-renormalization of $K_{1}$.
Then ${\cal K} = \{ K_{i}\}_{i=1}^{\infty}$ is a sequence of
renormalizations of $K_{1}$. Let $\{ (U_{i},V_{i}, f_{i})\}_{i=1}^{\infty}$
be a sequence of renormalizations with filled-in Julia set $K_{i}$ where
$f_{i} =f_{i-1}^{\circ n_{i-1}}$ and where $n_{i-1}$ is the period of
the critical tableau $T^{i-1}(0)$ of the two-dimensional Yoccoz puzzle
for $(U_{i-1},
V_{i-1}, f_{i-1})$, for $i\geq 2$. Suppose $U_{i}$ is a renormalization of
$(U_{i-1}, V_{i-1}, f_{i-1})$. We describe $(U_{1}, V_{1}, f_{1})$ as
$(n_{1}, n_{2}, \ldots)$-infinitely renormalizable. The grid $\{
T^{i}(0)\}_{i=1}^{\infty}$ is called
the three-dimensional critical tableau for $(U_{1},
V_{1}, f_{1})$.
Let $c(n)=f_{1}^{\circ n}(0)$. The critical orbit of $f_{1}$ is
$CO=\{
c(n)\}_{n=0}^{\infty}$. Let $GCO =\cup_{k=0}^{\infty}\cup_{n=0}^{\infty}
f_{1}^{-k}(c(n))$ be the grand critical orbit of $f_{1}$.

\vskip5pt
\proclaim Definition 2. An infinitely
renormalizable quadratic-like map
$(U_{1}, V_{1}, f_{1})$ has {\sl a priori} complex bounds if
there are a constant $\lambda >0$ and
a sequence of renormalizations $\{ (U_{i_{k}}, V_{i_{k}},
f_{i_{k}})\}_{k=1}^{\infty}$ of $f_{1}$ such that
$$\hbox{mod}(V_{i_{k}}\setminus U_{i_{k}})\geq \lambda$$
for all $k\geq 1$. 

\vskip5pt
\proclaim Theorem 15. Suppose $(U_{1}, V_{1}, f_{1})$ is an
infinitely
renormalizable quadratic-like map having the {\sl a priori} complex 
bounds. Its filled-in Julia set $K_{1}$ is
locally connected at every point in $GCO$.

{\bf Proof.} Suppose, without loss of generality, that
$\{ U_{i}, V_{i}, f_{i}\}_{i=1}^{\infty}$
is the sequence of renormalizations in Definition 2. Let $\lambda >0$
be the constant in Definition 2.
Then $\{ U_{i}\}_{i=1}^{\infty}$
is a sequence of nested domains containing $0$.
Consider the annulus $A_{i} =\overline{U}_{i}\setminus U_{i+1}$
for $i\geq 1$. For each $i\geq 1$, let $cv_{i} =f_{i}(0)$ and let $\gamma_{i}$
be a curve in $V_{i}\setminus f_{i}^{\circ 2}(U_{i+1})$ connecting $cv_{i}$ and a point on the boundary
of $V_{i}$. Let $0\in U_{i}'\subset U_{i}$ be the connected component of
the pre-image of $V_{i}\setminus \gamma_{i}$ under $f_{i}^{\circ 2}$.
Then $f_{i}^{\circ 2}:
U_{i}'\rightarrow V_{i}\setminus \gamma_{i}$ is a degree two
branch covering.
Moreover, $f_{i}:U_{i}' \rightarrow f_{i}(U_{i}')\subset U_{i}$ is
also a degree two branch covering map. Thus $f_{i}: U_{i}\setminus
U_{i}' \rightarrow
V_{i}\setminus f(U_{i}')$ is a degree two branch covering map. This implies
that $$ \hbox{mod}(U_{i}\setminus U_{i}') = {1\over 2}\cdot
\hbox{mod}(V_{i}\setminus f_{i}(U_{i}')).$$
But $V_{i}\setminus U_{i}$ is a sub-annulus of $V_{i}\setminus f_{i}(U_{i}')$.
So
$$ \hbox{mod}(U_{i}\setminus U_{i}') \geq {1\over 2}\cdot
\hbox{mod}(V_{i}\setminus U_{i}) > {\lambda \over 2}.$$
Remember that $U_{i+1}$ is the domain of the renormalization
$f_{i+1}=f_{i}^{n_{i}}$ where $n_{i}\geq 2$. We have $U_{i+1}\subset U_{i}'$.
Hence $U_{i}\setminus U_{i}'$ is a sub-annulus $A_{i}$. So
$$\hbox{mod}(A_{i})> {\lambda\over 2}.$$
Let $A_{\infty} =\cap_{i=1}^{\infty} U_{i}$.
Since $U_{1} \setminus A_{\infty} =\cup_{i=1}^{\infty} A_{i}$,
$$ \hbox{mod}(U_{1}\setminus A_{\infty})\geq  \sum_{i=1}^{\infty}
\hbox{mod}(A_{i})
=\infty.$$ Thus, $A_{\infty} =\{ 0 \}$. This implies that the diameter
$d(U_{i})$ tends to $0$ as $i$ goes to infinity.

Let $\alpha_{i}$ be the separate fixed point of $f_{i}$. Let
$\Gamma_{i}$ be a cycle of periodic external rays of $f_{1}$
landing at $\alpha_{i}$ (refer to the end of \S 3). Let $\eta^{i}_{0}$
be the set consisting of the
closures of the connected components
of $U_{1}\setminus \Gamma_{i}$ and let $\eta^{i}_{n} =f_{1}^{-n}(\eta^{i}_{1})$ for
$n\geq 1$. Then $\{ \eta^{i}_{n} \}_{n=0}^{\infty}$ is a puzzle for
$f_{1}$. Let $C^{i}_{n}$ be the member of $\eta^{i}_{n}$ containing $0$.
Consider the critical end
$$0\in \cdots \subseteq C_{n}^{i} \subseteq C_{n-1}^{i} \subseteq \cdots
\subseteq C_{1}^{i}\subseteq C_{0}^{i}$$
and the critical tableau $T^{i}(0)=(a_{nm}^{0}(i))_{n\geq 0, m\geq 0}$.
Since $(U_{1}, V_{1}, f_{1})$ is $(n_{1}, n_{2}, \ldots
)$-renormalizable, $T^{i}(0)$ is periodic of
period $m_{i}=\prod_{j=1}^{i}n_{j}$. Thus there is an integer $N>0$ such that
$a_{nm}^{0}(i)=0$ for $n\geq N$ and $0< i< m_{i}$. Thus
$$f_{1}^{\circ m_{i}} = f_{i}^{\circ n_{i}}: \CC_{m_{i}+N}^{i}
\rightarrow
\CC_{N}^{i}$$ is a degree two proper holomorphic map. We may assume that
$C_{m_{i}+N}^{i} \subset \CC_{N}^{i}$ (otherwise, we can modify
$C_{n}^{i}$ as the proof of Theorem 9). Therefore,
$$f_{i+1}=f_{1}^{\circ m_{i}}=f_{i}^{\circ n_{i}}: \CC_{m_{i}+N}^{i}
\rightarrow \CC_{N}^{i}$$
is a quadratic-like map. Since
$$K_{i+1} =\cap_{j=0}^{\infty} f_{1}^{-jm_{i}} (C^{i}_{N})=\cap_{j=0}^{\infty}
C^{i}_{jm_{i}+N},$$
there is a $C^{i}_{k(i)}$ contained in $U_{i}$. The diameter
$d(C^{i}_{k(i)})$ of $C^{i}_{k(i)}$ tends to zero as $i$ goes to infinity.
From Lemma 2, $C_{k(i)}^{i}\cap K_{1}$ is connected. So $\{
C_{k(i)}^{i}\}_{i=1}^{\infty}$ is a basis of connected neighborhoods at
$0$.

For any $x=f^{\circ n}_{1}(0)$, consider $\{ D_{i,n}(x) = f^{\circ
n}_{1}(C_{k(i)}^{i})\}_{i=1}^{\infty}$. It is a basis of connected
neighborhoods at $x$.
For any $y$ in $f^{-m}_{1}(x)$ (where $0$ is not in $f^{-n}(x)$ for
$0<n \leq m$), there is an open
neighborhood $W$ of $y$ such that $f_{1}^{\circ m}: W\rightarrow f^{\circ
m}(W)$ is a homeomorphism. Let $g$ be its inverse. Then $\{
g(D_{i,n}(x))\}_{i=1}^{\infty}$ is a basis of connected neighborhoods at $y$.
\hfill \bull

\vskip5pt
Suppose $(U_{1}, V_{1}, f_{1})$ is an $(n_{1}, n_{2}, \ldots
)$-infinitely
renormalizable quadratic-like map. We call the puzzle $\{ \{
\eta_{n}^{i}\}_{n=0}^{\infty} \}_{i=1}^{\infty}$ constructed in Theorem 15,
the three-dimensional Yoccoz puzzle for $(U_{1}, V_{1}, f_{1})$.
Let $m_{i}=\prod_{j=1}^{i} n_{j}$ and $\{ K_{i}\}_{i=1}^{\infty}$ be the
sequence of renormalizations of $K_{1}$. Let $c(j)=f^{\circ j}(0)$ for
$j\geq 0$, and let
$CO=\{ c(j)\}_{j=0}^{\infty}$ be the critical
orbit of $f_{1}$.

\vskip5pt
\proclaim Definition 3. An infinitely renormalizable quadratic-like
map $(U_{1}, V_{1}, f_{1})$ is unbranched if
there are a constant $\lambda >0$ and
a sequence of domains $\{ W_{k}\}_{k=1}^{\infty}$ such that
$W_{k}\supset K_{i_{k}}$, such that
$$\hbox{mod}(W_{k}\setminus K_{i_{k}})\geq \lambda,$$
and such that $W_{k}\cap CO= \{ c(jm_{i_{k}})\}_{j=0}^{\infty}$ for all $k\geq 1$

\vskip5pt
\proclaim Theorem 16. Suppose $(U_{1}, V_{1}, f_{1})$ is an
infinitely renormalizable unbranched quad\-ratic-like map having
the {\sl a priori} complex bounds. Then its filled-in Julia set $K_{1}$
is locally connected.

\vskip5pt
{\bf Proof.} Suppose, without loss of generality,  that $k=i$ and $i_{k}=i$ in Definition
3, and that $\{ U_{i}, V_{i},
f_{i}\}_{i=1}^{\infty}$
is a sequence of renormalizations in Definition 2. Let $\lambda >0$ be
a constant satisfying Definitions 2 and 3.

Let $\{ \{ \eta_{n}^{i}\}_{n=0}^{\infty} \}_{i=1}^{\infty}$ be
the three-dimensional Yoccoz puzzle for $(U_{1}, V_{1}, f_{1})$.
Let $\{ C_{k(i)}^{i}\}_{i=1}^{\infty}$ be the basis of
connected neighborhoods constructed in
Theorem 15. By choosing $k(i)$ large enough, we can have
$$\hbox{mod}(W_{i}\setminus C_{k(i)}^{i})
\geq {\lambda\over 2}$$
for all $i\geq 1$.

If $x=0$, Theorem 15 says that $K_{1}$ is locally connected at $x$.
For each $x\neq 0$ in $K_{1}$, there are two cases:
either (1) $x$ is non-recurrent, which means that there is an integer $i\geq 1$ such that
$\{ f_{1}^{\circ n}(x)\}_{n=0}^{\infty}\cap C_{k(i)}^{i}=\emptyset$;
or else (2) $x$ is recurrent.

In case (1), we prove it
by applying the results in hyperbolic geometry.
Let $x$ be a non-recurrent point. Then there is a $C^{i}_{k(i)}$
such that the orbit
$O(x)= \{ f_{1}^{\circ n} (x)
\}_{n=0}^{\infty}$ is disjoint from the interior of $C^{i}_{k(i)}$.
Let $r$ be an external ray of $f_{1}$ cutting $U_{1}\setminus
f_{1}(C^{i}_{k(i)})$ into a
simply connected domain $\Omega$.
Consider two inverse branches $Q_{1}$ and $Q_{2}$ of $f_{1}|\Omega$.

Let $\tilde{\Omega} = V_{1}\setminus (f_{1}(C^{j}_{k(j)})\cup r)$ for some $j>i$
such that $\overline{\Omega}\subset \tilde{\Omega}$. Let $\tilde{Q}_{1}$ and
$\tilde{Q}_{2}$ be the inverses of $f_{1}|\tilde{\Omega}$ and let
$\tilde{\Omega}_{1}$ and $\tilde{\Omega}_{2}$ be the images of
$\tilde{\Omega}$ under $\tilde{Q}_{1}$ and $\tilde{Q}_{2}$, respectively.
Consider $\tilde{\Omega}$, $\tilde{\Omega}_{1}$, and $\tilde{\Omega}_{2}$
as hyperbolic Riemann surfaces. Then $Q_{1}: \Omega \rightarrow \Omega_{1}$
and $Q_{2}: \Omega \rightarrow \Omega_{2}$
strictly contract the hyperbolic distances of
$\tilde{\Omega}$, $\tilde{\Omega}_{1}$, and $\tilde{\Omega}_{2}$.
Since $\Omega$ contains finite number of connected components
of $K_{1} \setminus f_{1}(C^{i}_{k(i)})$, we can cut
$\Omega$ into finite number of
simply connected domains $\Omega_{1}'$, $\ldots$, $\Omega_{l}'$ such
that
$\Omega_{k}'\cap K_{1}$ is connected for each $1\leq k\leq l$. Since the
orbit $O(x)$ is contained in $\Omega$,
from the images of $\Omega_{1}'$, $\ldots$, $\Omega_{l}'$ under
the semi-group generated by $Q_{1}$ and $Q_{2}$,
we can get a basis of connected neighborhoods at $x$.
Therefore, $K_{1}$ is locally connected at $x$.

In case (2), $f_{1}^{\circ n}(x)$ enters $C_{k(i)}^{i}$ infinitely many
times. For each $i\geq 1$,
consider the puzzle $\eta^{i}=\{ \eta_{n}^{i}\}_{n=0}^{\infty}$. For the $x$-end,
$$ x\in \cdots \subseteq D^{i}_{n}(x) \subseteq D_{n-1}^{i}(x) \subseteq
\cdots \subseteq D^{i}_{1}(x) \subseteq D^{i}_{0}(x),$$
let $T^{i}(x)=(a_{nm}(i))_{n\geq
0, m\geq 0}$ be the corresponding
tableau. Let $q_{i}$ be the integer such that
$a_{k(i)q}=0$ for $0\leq q< q_{i}$ and such that $a_{k(i)q_{i}}=1$. Let
$p_{i}=k(i)+q_{i}$.
Then $g_{i,x}=f_{1}^{\circ q_{i}}: D^{i}_{p_{i}}(x)
\rightarrow C_{k(i)}^{i}$ is a proper
holomorphic diffeomorphism. Since there is no critical point in
$W_{i}\setminus C_{k(i)}^{i}$, this diffeomorphism can be extended to a
proper holomorphic diffeomorphism on $W_{i}$. Let
$x_{i}=f_{1}^{\circ q_{i}}(x)$.
We can modify $W_{j}$ such that
$${\lambda \over 2} \leq \hbox{mod}(W_{j}\setminus C_{k(j)}^{j}) \leq  2\lambda.$$
From the previous theorem, the diameter
$d(C_{k(j)}^{j})$ tends to zero as $j$ goes to infinity. We can find an integer
$j>i$ such that $W_{j}\subset C_{k(i)}^{i}$. Let
$$x_{j}=f_{1}^{\circ q_{j}}(x)=f_{1}^{\circ (q_{j}-q_{i})}(x_{i}).$$
Consider the puzzle $\eta^{j}=\{ \eta_{n}^{j}\}_{j=1}^{\infty}$. Let
$$ x_{i} \in \cdots \subseteq D^{j}_{n}(x_{i}) \subseteq
D_{n-1}^{j}(x_{i}) \subseteq \cdots \subseteq
D^{j}_{1}(x_{i}) \subseteq D^{j}_{0}(x_{i})$$
be the $x_{i}$-end and let
$T^{j}(x_{i})=(b_{nm}(j))_{n\geq
0,m\geq 0}$ be the corresponding tableau. Then one can check that
$f_{1}^{\circ (q_{j}-q_{i})}$ is a proper holomorphic diffeomorphism from
$$D^{j}_{k(j)+q_{j}-q_{i}}(x_{i})\rightarrow C_{k(j)}^{j}.$$ In other
words,
$b_{k(j)m}=0$ for $0\leq m< q_{j}-q_{i}$, but $b_{k(j)(q_{j}-q_{i})} =1$.
Let $g_{ij}$ be the inverse of $f_{1}^{\circ (q_{j}-q_{i})}:
D^{j}_{k(j)+q_{j}-q_{i}}(x_{i})\rightarrow C_{k(j)}^{j}$.
Then $g_{ij}$ can be extended to $W_{j}$. Since
$C_{k(i)}^{i}$ is bounded by external rays landing at some
pre-iamges of $\alpha_{i}$ under iterations of $P_{c}$ and
is a part of an invariant net under $f_{1}$, $W_{ij}=g_{ij}(W_{j})$ is
contained in $C_{k(i)}^{i}$. Thus
$$\hbox{mod}(W_{i}\setminus W_{ij}) \geq {\lambda \over 2}.$$
Consider $X_{i}=g_{i,x}(W_{i})$ and
$X_{j}=g_{j,x}(W_{j})=g_{i,x}(W_{ij})$.
Then
$$\hbox{mod} (X_{i}\setminus
X_{j}) \geq {\lambda\over 2}.$$
Therefore, inductively, we can find a sequence of
domains $\{ X_{i_{s}}\}_{s=1}^{\infty}$ such that
$$\hbox{mod}(X_{i_{s}}\setminus X_{i_{s+1}}) \geq {\lambda\over 2}$$
for $s\geq 1$. Thus the diameter of $X_{i_{s}}$ tends to zero as $s$ goes to
infinity.
For each puzzle $\eta^{i_{s}}=\{ \eta_{n}^{i_{s}}\}_{n=0}^{\infty}$,
consider both the $x$-end,
$$ x\in \cdots \subseteq D^{i_{s}}_{n}(x) \subseteq D_{n-1}^{i_{s}}(x)
\subseteq \cdots \subseteq
D^{i_{s}}_{1}(x) \subseteq D^{i_{s}}_{0}(x)$$
and the corresponding
tableau $T^{i_{s}}(x)=(a_{nm}(i_{s}))_{n\geq 0, m\geq 0}$.
For each $k(i_{s})$, there is an integer $q_{i_{s}}$ such that
$a^{i}_{k(i_{s})q} =0$
for $0 \leq q< q_{i_{s}}$ and such that $a^{i}_{k(i_{s})q_{i_{s}}}=1$. Then for
$p_{i_{s}} =k(i_{s})
+q_{i_{s}}$, $f_{1}^{\circ q_{i_{s}}}: D_{p_{i_{s}}}(x) \rightarrow
C_{k(i_{s})}$ is a proper holomorphic diffeomorphism. This implies that
$$D_{p_{i_{s}}}^{i_{s}}(x) \subseteq X_{i_{s}}.$$
So the diameter $d(D_{p_{i_{s}}}^{i_{s}}(x))$ tends to zero as $s$ goes to infinity.
From Lemma 2, $\{ D_{p_{i_{s}}}^{i_{s}}(x)\}_{s=1}^{\infty}$ forms a
basis of connected neighborhoods at $x$.
\hfill \bull

\section {A Generalized Sullivan's Sector Theorem}

Let $I=[0,1]$ be the closed unit interval. Let ${\cal E}_{0}$ be the set of all
functions $G$ such that

\vskip3pt
{\bf (1)} $G:I_{G} \supseteq I \rightarrow G(I_{G})$ is a homeomorphism
and $G(0)=0$, $G(1)=1$, and

\vskip3pt
{\bf (2)} $G$ can be extended to be a schlicht function $g$ on
$\bold C_{G}=(\bold C \setminus \bold R^{1})\cup
\II_{G}$ preserving upper- and lower-half planes.

\vskip3pt
\noindent We assume $I_{G}$ is the maximum
interval satisfies {\bf (1)} and {\bf (2)} for each $G$ in ${\cal E}_{0}$
and call it the definition interval of $G$. We will not distinguish $g$ and
$G$ anymore.
Take
$$S_{\gamma}(z) = r^{1\over \gamma} e^{{\theta \over \gamma}i}:
\bold C\setminus \{ x<0\}\rightarrow \bold C$$
as the standard $\gamma$-root
where $z=re^{\theta i}$ for $r>0$, $-\pi < \theta < \pi$, and
$\gamma >1$.
For every $a \leq 0$, we call $L_{a}(z) =E S_{\gamma}(z-a)+F$ a
$\gamma$-root at $a$ where
$E=1/ \Big( (1-a)^{1\over \gamma}-(-a)^{1\over \gamma}\Big)$ and $F=-(-a)^{1\over \gamma}E$
are determined by $L_{a}(0) =0$ and $L_{a}(1)=1$.
Then $L_{a}$ is an element in ${\cal E}_{0}$ whose definition interval is $[a, \infty)$.

Suppose that $L_{a}$ is a $\gamma$-root at $a$ and that $G$ is an
element in ${\cal
E}_{0}$ whose definition interval is $I_{G}$. We say that $L_{a}$ and
$G$ are compatible
if $[a, 1] \subset G(I_{G})$. For a compatible pair $L_{a}$ and $G$,
let $a'=G^{-1}(a)$,
$J=[a',1]$, $L\cup R= I_{G}\setminus J$, and $b=\min \{ |L|, |R|\}$.

\vskip5pt
\centerline{\psfig{figure=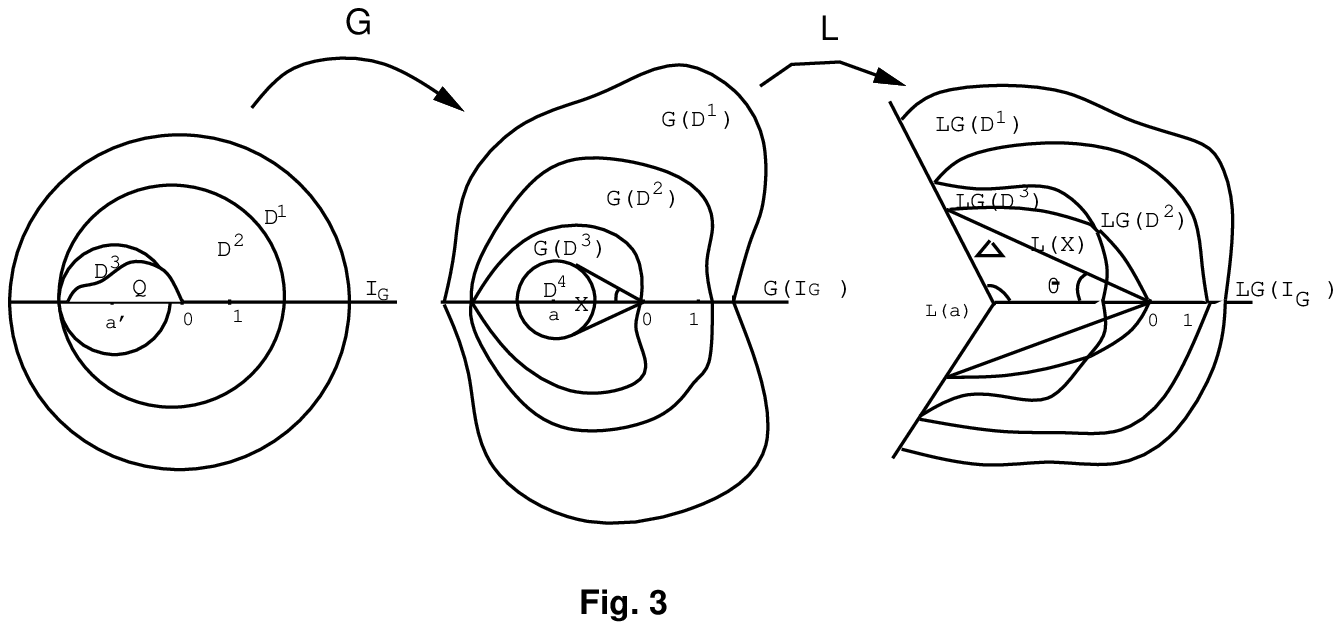}}
\vskip5pt

We consider several disks related to
a compatible pair
$L_{a}$ and $G$ (see Fig. 3).
Set $D^{1}$
to be the closed disk centered at
the middle point $(1+a')/2$ of $[a',1]$ with diameter $1+|a'|+2b$. Set
$D^{2}$ to be the closed disk centered at
the middle point $(1+a')/2$ of $[a',1]$ with diameter $k=\min\{
1+|a'|+b, 2(1+|a'|)\}$.
Set
$D^{3}$ to be the maximum closed
disk centered at $a'$ and contained in $D^{2}$. Then the radius
of $D^{3}$ is $d=\min \{ b/2, (1+|a'|)/2\}$.
We note that
$$D^{3} \subset D^{2} \subset D^{1}.$$

The map $G$ is a schlicht function on $D^{1}$.
Let $\mu =(1+|a'|)/b$
and $\nu=(2\mu +3)^{4}$.
From Koebe's distortion theorem (see~[BIE])
for any $\xi$ and $\eta$ in $D^{2}$,
$$ \nu^{-1} \leq {|G'(\xi)| \over |G'(\eta)|} \leq \nu.$$
Using the fact that $G(I)=I$, there is at least one $\eta$ in $D^{2}$ such
that $|G'(\eta)|=1$. Hence
$$ \nu^{-1} \leq |G'(\xi)| \leq \nu$$
for all $\xi$ in $D^{2}$.
This implies that
$$ {|a'|\over \nu} \leq |a| \leq |a'|\cdot \nu .$$

Let ${\bf UH} =\{ z=x+yi \in
\bold C\; |\; y >0\}$ be the upper-half plane. For any
$z$ in ${\bf UH}$, let $\theta (z) = \arg (z)$.

\vskip5pt
\proclaim Theorem 17.
Suppose that $L=L_{a}$ and $G$ are a pair of compatible elements in
${\cal E}_{0}$.
There is a constant $0< \theta < \pi $
depending only on $\mu$
such that
the image $L(G({\bf UH}))$ of the upper-half plane under $L\circ G$
contains an open triangle $\Delta$
based on $[L(a),0]$ whose angle at $L(a)$ is $\pi/\gamma$ and whose angle
at $0$ is $\theta$.

\vskip5pt
{\bf Proof.}
The image $G(D^{3})$ of $D^{3}$ under $G$
contains the
closed disk $D^{4}$ centered at $a$ with radius $d/\nu$. Similarly, for any
$a'\leq x\leq 1$, consider the closed disk $D(x)$ centered at $x$ with
radius $d$. Then $D(x) \subset D^{2}$ and $G(D(x))$ contains the closed disk
centered at $G(x)$ with radius $d/\nu$. This implies that the convex-hull
$X$ of $\{ 0\}\cup D^{4}$ is contained in $G(D^{2})$.
Either $X=D^{4}$ or $X\cap {\bf UH}$ has an angle at $0$. In the later case,
the angle $\varphi$ of $X\cap {\bf UH}$ at $0$ has $\sin \varphi
=d/(|a| \nu)$.

Since $L$ is a $\gamma$-root for $\gamma >1$, the convex set
$L(X)$ contains a triangle $\Lambda$ based on $[L(a), 0]$
whose angle at $L(a)$ is $\pi/\gamma$ and whose angle $\omega$
at $0$ can be calculated from
$${\sin \omega \over \sin ({\pi \over \gamma}+\omega)} =
\Big({d\over |a| \nu} \Big)^{1\over \gamma}$$
through the law of sines. Because $d=\min\{ b/2, (1+|a'|)/2\}$,
$${d\over |a|} \geq {d\over |a'| \nu} \geq \min \{ {1\over
2 \mu \nu}, {1\over 2\nu}\}.$$
Hence $\Lambda$ contains a triangle $\Delta$
based on $[L(a), 0]$ whose angle at $L(a)$ is still $\pi/\gamma$
and whose angle $\theta$ at $0$ is calculated from
$${\sin \theta \over \sin ({\pi \over \gamma}+\theta)} =
\Big( \min \{ {1\over 2 \mu \nu},
{1\over 2\nu} \} \Big)^{1\over \gamma}.\eqno
\bull
$$

\vskip5pt
Suppose $Q=(L\circ G)^{-1} (\overline{\Delta})$, where $\Delta$
is the triangle obtained in Theorem 17. Then
$$Q \subset D^{2} \subset D^{1}.$$

Suppose $\{ (L_{i}, G_{i})\}_{i=0}^{n}$
is a sequences of compatible pairs
in ${\cal E}_{0}$ where $L_{i}$ is a $\gamma$-root at $a_{i}$.
Let $a_{i}'$, $b_{i}$, $k_{i}$, $d_{i}$, $\mu_{i}$, $\nu_{i}$, and $\theta_{i}$
be the numbers, and $D_{i}^{1}$, $D_{i}^{2}$, $D_{i}^{3}$, $\Delta_{i}$, and $Q_{i}$ be the sets
corresponding to each compatible pair $L_{i}$ and $G_{i}$.
Let
$$ {\cal L} = L_{n}\circ G_{n}\circ \cdots \circ L_{i}\circ G_{i}\circ
\cdots \circ L_{0}\circ G_{0}.$$
Then ${\cal L}$ is a schlicht function defined on
$\bold C_{I}=(\bold C\setminus \bold R^{1})\cup \II$.

\vskip5pt
\proclaim Definition 4. We call ${\cal L}$ a root-like
map if there are constants $C>0$ and $\lambda >1$ such that

\vskip3pt
{\bf (i)} $a_{0}=0$ and $a_{1} \geq 1/C$,

\vskip3pt
{\bf (ii)} $|a_{j}| \geq
\max \{ (\lambda^{j-i}/C) \cdot |a_{i}|, \Big( 1+(\lambda
-1)/C\Big)\cdot |a_{i}|\} $ for all $1\leq i<j \leq  n$, and

\vskip3pt
{\bf (iii)} $\mu_{i} <C$ for all $1\leq i\leq n$.

\vskip5pt
\proclaim Theorem 18. Suppose ${\cal L}$ is a
root-like map. There is a constant $\theta>0$ depending only on
$\lambda$ and $C$ such that the image of the upper-half plane
under ${\cal L}$ is contained in the
sector
$${\bf Sec}_{\theta}=\{ z \in \bold C\;|\; 0 \leq \arg (z) \leq \pi
-\theta \}.$$

\vskip5pt
Before we prove this theorem, we introduce some basic
results in hyperbolic geometry. Let
$\bold C_{I}=(\bold C\setminus \bold R^{1})\cup \II$ be a plane domain.
Then $q(z) =-z^{2} /(1-z^{2})$ is
a diffeomorphism from the upper-half plane ${\bf UH}$ onto $\bold C_{I}$.
Consider ${\bf UH}$ to be a hyperbolic plane with Poincar\'e metric
$d_{H, {\bf UH}}s=|dz|/y$
for $z=x+yi$. This metric induces a hyperbolic metric
$d_{H, \bold C_{I}}s = q_{*} (|dz|/y)$ on ${\bf
C}_{I}$.
The plane domain $\bold C_{I}$ under this metric
is a hyperbolic Riemann surface.
Let $d=d_{H, \bold C_{I}}$ be the induced hyperbolic distance.
We note that $q$ maps the
positive imaginary line in ${\bf UH}$ onto the interval $I$ and maps the
real line, which is the boundary of ${\bf UH}$, onto the set ${\bf
R}^{1} \setminus \II$.

\vskip5pt
\proclaim Lemma 6.
A hyperbolic neighborhood $\Phi (r)=\{z \in \bold C_{I} \; |\; d(z,
I)<r\}$ is the union of two Euclidean disks $D^{+}$
and $D^{-}$, symmetric to each other with
respect to $I$, centered at $c^{+}$ and at $c^{-}=-c^{+}$,
with the same radius $R^{+}=R^{-}$. Moreover,
$$ R^{+}= {1\over 2 \sin \beta} \qquad \hbox{and} \qquad
c^{+} = {1\over 2}+ {\cot \beta  \over 2}i $$
where
$$\beta = 4 \cot^{-1} (e^{r})$$
is the angle at $0$ between $\partial D^{+}$
and the negative real line (and the angle at $1$ between $\partial
D^{+}$ and the ray $[1, \infty )$).

\vskip5pt
{\bf Proof.} Consider the pre-image $\Phi'= q^{-1}(\Phi (r))$. It is
a hyperbolic neighborhood in ${\bf UH}$
and consists of all points in ${\bf UH}$ whose hyperbolic distances
to the half-line $l_{+}=\{ z=yi\; | \; y>0 \}$ are less
than $r$. The boundary $\partial \Phi'$ consists of two rays starting
from
$0$. Thus $\Phi'$ is a sector, symmetric with respect to the half-line
$l_{+}$.
Suppose $\beta/2$ is the outer angle of this sector (with
respect to the real line). Since a geodesic in
${\bf UH}$ is a semi-circle or half-line perpendicular to the real line,
it is easy to check that
$$\log \big( \cot {\beta\over 4}\big) = r.$$
Therefore,
$\Phi (r)$ is the union of two disks $D^{+}$ and
$D^{-}$ symmetric with respect
to $I$. The angle between $\partial D^{+}$ (or $\partial D^{-}$) and the
negative real line is $\beta$.
Moreover, every point $z$ in $\partial D^{+}$ (or
$\partial D^{-}$) views $I$ under the same angle $\beta$, that is, every
triangle $\triangle(0z1)$ has the angle $\beta$ at $z$. Now consider the point
$u$ such that the segment $\overline{1u}$ is a diameter
of $D^{+}$. The triangle $\triangle(u01)$ is a
right triangle. We can calculate the length $2R^{+}$ of the segment
$1u$ and length $|u|$ of the segment $\overline{0u}$ as follows:
$$ 2 R^{+} = {1 \over \sin \beta}, \hskip40pt |u| = \cot \beta.$$
Therefore,
$$ R^{+}= {1\over 2 \sin \beta} \qquad \hbox{and} \qquad c^{+}
= {1\over 2}+ {\cot \beta \over 2}i .\eqno
\bull
$$

\vskip5pt
\proclaim Lemma 7. Let $z$ be a point in $\bold C_{I}$ and let
$\Phi(r)=D^{+}\cup D^{-}$ be the smallest hyperbolic neighborhood
containing $z$. The Euclidean radius $R^{+}$ of $D^{+}$ (and
$D^{-}$) is
$${|z-1|\over 2 \sin (\arg (z))}.$$

\vskip5pt
{\bf Proof.} Let $u$ be the point in $\partial D^{+}$ such that the
segment $\overline{1u}$ is a diameter of $D^{+}$.
The angle of the triangle $\triangle(0z1)$ at $z$ and the angle of the
triangle $\triangle(0u1)$
at
$u$ are both $\beta$ (see the previous lemma). Now applying the law of sines,
$$ {\sin (\arg (z)) \over |1-z|} = \sin \beta = {1\over |1-u|}.$$
Therefore,
$$2 R^{+}= {|z-1|\over \sin (\arg (z))}.
\eqno \bull$$

\vskip5pt
Suppose ${\cal L}$
is a root-like map.
Then there is a constant angle $\sigma$ and a constant $C_{0}>0$ depending
only on $C$ such
that
$\theta_{i} \geq \sigma$ and $\nu_{i} \leq C_{0}$ for all $0\leq i \leq
n$.

For each $1\leq i\leq n-1$, let
$A_{i+1} = Q_{i+1}\setminus \Delta_{i}$. Let $A_{1}=Q_{1}
\setminus \{ z|\in \bold C, 0\leq \arg (z) < \pi - \pi/\gamma\}$.
Suppose $\Phi (r_{i})=D^{+}_{i}\cup D^{-}_{i}$ is the smallest
hyperbolic neighborhood
in $\bold C_{I}$ containing $A_{i}\not= \emptyset$ and let $D^{0}_{i}$
be the smallest disk centered at $1/2$ containing $\Phi (r_{i})$.
Let $R^{+}_{i}$ be the Euclidean radius of $D^{+}_{i}$ and let $R_{i}$
be the Euclidean radius of $D^{0}_{i}$.
Since $z$ and $1$ are in $D_{i}^{2}$, from Lemmas 6 and 7,
$$ R_{i} \leq 2\cdot R^{+}_{i} \leq {k_{i} \over \sin \theta_{i}} 
\leq {2\cdot (1+|a_{i}'|) \over \sin \theta_{i}} \leq
{2\cdot (1+C)\cdot |a_{i}'| \over \sin \sigma} =C_{1}\cdot |a_{i}'|$$
where $k_{i}= \min\{ 1+|a_{i}'|+b_{i}, 2(1+|a_{i}'|)\}$ is the diameter of
$D_{i}^{2}$ and where $C_{1}= 2\cdot (1+C)/\sin \sigma$ is a positive constant depending
only on $C$ (see Fig. 4).

For every $1\leq i \leq n$, the definition interval of
$\phi_{i}= L_{i}\circ G_{i}$ is $[a_{i}', a_{i}'']= [a_{i}', \infty)\cap I_{G_{i}}$. According to
${\bf (iii)}$ of Definition 4, the right-endpoint
$a_{i}''$ satisfies
$$a_{i}'' \geq {1+|a_{i}'|\over C}
+1 \geq 1+ {|a_{i}'|\over C}.$$
For every $1\leq i<j \leq n$,
$$ |a_{j}'| \geq {\lambda^{j-i} \over C C_{0}^{2}} |a_{i}'|.$$
Therefore, if $\tau = \min\{ 1, 1/ (C\cdot C_{0})^{2}\}$, we have
$$[a_{j}', a_{j}''] \supset I_{i,\tau}=[\lambda^{j-i}\tau a_{i}',
\lambda^{j-i} \tau |a_{i}'|+1].$$
In other words, $\phi_{j}=L_{j}\circ G_{j}$ are schlicht functions on
$\bold C_{I_{i,\tau}}=(\bold C\setminus \bold R^{1})\cup I_{i, \tau}$
for all $1\leq i< j \leq n$.

\vskip5pt
\centerline{\psfig{figure=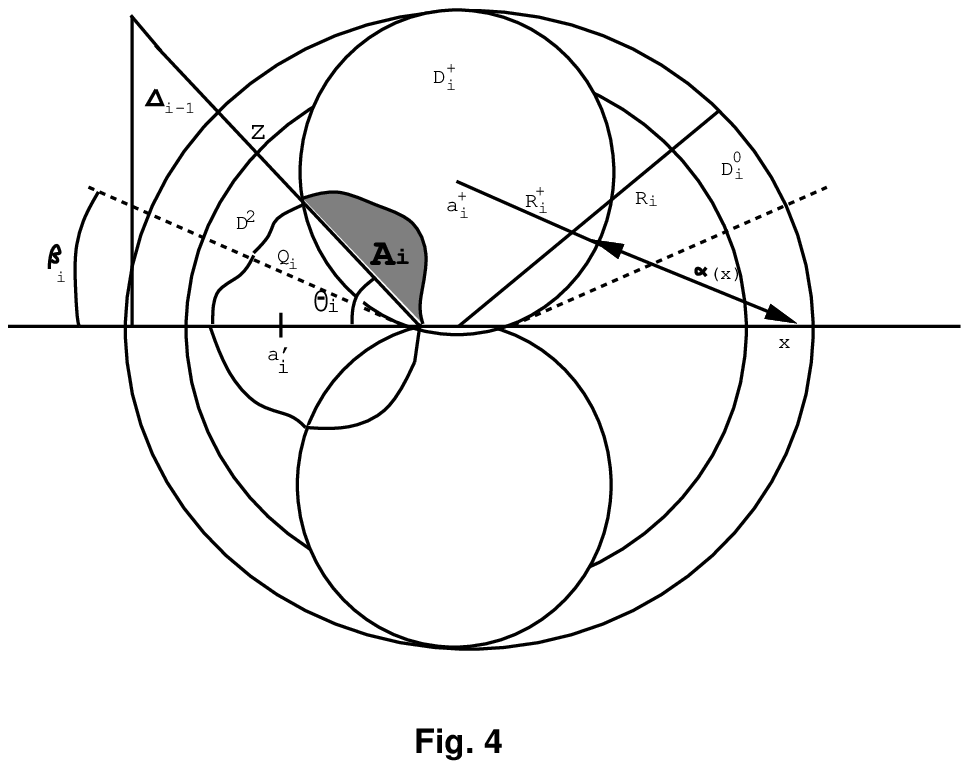}}
\vskip5pt

\proclaim Lemma 8. There is a fixed integer $n_{0}>0$ depending only on $C$
and on $\lambda$ such that for any
$0\leq i < n$ and any $j \geq i+ n_{0}$,
$$ |a_{j}'|+{1\over 2} >\lambda^{j-i-n_{0}} R_{i} \qquad \hbox{and} \qquad
a_{j}''-{1\over 2} > \lambda^{j-i-n_{0}} R_{i}.$$

\vskip5pt
{\bf Proof.}
Take $m_{0}$ as the biggest integer such that
$$ m\leq { \log ({C_{1}\over \tau}) \over \log \lambda
}.$$
From the above estimates, one can see that $n_{0} =\max\{m_{0}+1, 0\}$ is the
integer satisfying the lemma.
\hfill \bull

\vskip5pt
\proclaim Lemma 9. There is a constant $C_{2}>0$ depending only on $\lambda$ such that for any
$0\leq i < n-n_{0}$, let $\Pi_{i} =L_{n}\circ G_{n} \circ \cdots
\circ L_{j}\circ G_{j} \circ \cdots \circ L_{i+n_{0}+1}\circ
G_{i+n_{0}+1}$,
the distortion of $\Pi_{i}$ on $\Phi (r_{i})$ is
bounded by $C_{2}$, more precisely,
$$ {|\Pi_{i}' (\xi)| \over |\Pi_{i}' (\eta)|} \leq C_{2}$$
for all $\xi$ and $\eta$ in $\Phi (r_{i})$.

\vskip5pt
{\bf Proof.} Consider the disk $D^{5}_{j}$ centered at $1/2$ with radius
$t_{j}=\min\{|a_{j}'|+1/2, a_{j}''-1/2\}$. The composition $\phi_{j}=L_{j}\circ G_{j}$ is a schlicht
function on $D^{5}_{j}$. Let $r_{ij}=R_{i}/t_{j}$ be the
ratio of the radii of $D_{i}^{0}$ and $D_{j}^{5}$. Then
$$r_{ij} \leq \lambda^{i+n_{0}-j}$$
for $i+n_{0}+1 \leq j \leq n$. From
Koebe's distortion theorem (see~[BIE]),
$${|\phi_{j}' (\xi)| \over |\phi_{j}' (\eta)|} \leq \Big( {1+r_{ij}\over
1-r_{ij}}\Big)^{4} \leq \Big( {1+\lambda^{-k} \over 1-\lambda^{-k}}
\Big)^{4}$$
for $0< k=j-i-n_{0}\leq n-i-n_{0}$ and for all $\xi$ and $\eta$ in
$D^{0}_{i}$.
Since $\phi_{j}$ is a schlicht function on $\bold C_{I}$, it contracts
the hyperbolic distance $d$ on $\bold C_{I}$. Thus
$$\phi_{j} (\Phi (r_{i})) \subset
\Phi (r_{i}) \subset D^{0}_{i}.$$
Therefore, by the chain rule,
$${|\Pi_{i}' (\xi)| \over |\Pi_{i}' (\eta)| } = \prod_{j=n_{0}+1}^{n}
{|\phi_{j}'(\xi)|\over |\phi_{j}'(\eta)|} \leq
C_{2} = \Big( \prod_{k=1}^{\infty} {1+\lambda^{-k} \over 1-\lambda^{-k}}
\Big)^{4}$$
for all $\xi$ and $\eta$ in $\Phi(r_{i})$.
\hfill \bull

\vskip5pt
\proclaim Lemma 10. There is a constant $C_{3}>0$ depending only on $C$ such that for
$0\leq i\leq n$, for $i\leq j\leq i+n_{0}$, for $\phi_{j} = L_{j}\circ
G_{j}$, and for all $\xi$ and $\eta$ in $D^{+}_{i}$ (or $D^{-}_{i}$),
$$ {|\phi_{j}'(\xi)| \over |\phi_{j}' (\eta)|} \leq C_{3}.$$

\vskip5pt
{\bf Proof.} Suppose $x$ is a real number. Let
$\alpha=\alpha(x) =\min \{ |x-z|\; |\; z \in D^{+}_{i}\cap {\bf UH}\}$.
Suppose that $c=c^{+}_{i}$
and that $R=R^{+}_{i}$ are the center and the radius of $D^{+}_{i}$.
Suppose
$h=h^{+}_{i}$ is the length of the segment $\overline{c{1\over 2}}$ (the
straight line connecting $c$ and $1/2$). Since
the two triangles $\triangle(x{1\over 2}c)$ and $\triangle(0{1\over
2}c)$ are both
right triangles, then
$$ (\alpha + R)^{2}=({1\over 2}-x)^{2} +h^{2}$$
and
$$ R^{2} = ({1\over 2})^{2} + h^{2}.$$
Therefore,
$$ {\alpha + R \over R} =  \sqrt{1+{x^{2}-x \over R^{2}}}.$$
This implies that there is a constant $0< C_{4} < 1$ depending only on $C$ such
that for $x=a_{j}'$ or $x=a_{j}''$,
$$ {R\over \alpha + R} \leq C_{4}.$$
Now consider the largest disk $D^{6}_{j}$ centered at $c$ such that
$\phi_{j}$ is a schlicht function on it. Then the radius of $D^{6}_{j}$ is
greater than or equal to $\min \{ R+\alpha (a_{j}'), R+\alpha (a_{j}'')\}$.
From Koebe's distortion theorem (see~[BIE]),
$$ {|\phi_{j}' (\xi)| \over |\phi_{j}' (\eta)|} \leq C_{3}= \Big(
{1+C_{4}\over 1-C_{4}}\Big)^{4}$$
for any $\xi$ and $\eta$ in $D_{i}^{+}$.
\hfill \bull

\vskip5pt
Combining Lemmas 9 and 10, we obtained the following estimate:

\proclaim Lemma 11. There is a constant $C_{5} >0$ depending on $\lambda$ and
on $C$ such that for any
$0\leq i < n$, let $\Sigma_{i} =L_{n}\circ G_{n} \circ \cdots
\circ L_{j}\circ G_{j} \circ \cdots \circ L_{i}\circ
G_{i}$,
the distortion of $\Sigma_{i}$ on $D^{+}_{i}$ (or $D_{i}^{-}$) is
bounded by $C_{5}$, more precisely,
$$ {|\Sigma_{i}' (\xi)| \over |\Sigma_{i}' (\eta)|}\leq C_{5}$$
for all $\xi$ and $\eta$ in $D^{+}_{i}$ (or $D_{i}^{-}$).

\vskip5pt
{\bf Proof.} Since each $\phi_{j}$ is a schlicht function
on $\bold C_{I}$, it contracts the hyperbolic distance $d$ on ${\bf
C}_{I}$. So
$$ \phi_{j} (D^{+}_{i}) \subset D^{+}_{i}$$
for all $i\leq j \leq n$.
If $n-i\leq n_{0}$, then from Lemma 10 and the chain rule,
$$ {|\Sigma_{i}' (\xi)| \over |\Sigma_{i}' (\eta)|}\leq C_{3}^{n_{0}}$$
for all $\xi$ and $\eta$ in $D^{+}_{i}$.

Now we consider $n-i > n_{0}$ and write $\Sigma_{i} = \Pi_{i}
\circ \Theta_{i}$ where
$\Theta_{i} = L_{i+n_{0}}\circ G_{i+n_{0}}\circ \cdots \circ L_{i}\circ
G_{i}$ and where $\Pi_{i} = L_{n}\circ G_{n}\circ \cdots \circ
L_{i+n_{0}+1}\circ G_{i+n_{0}+1}$.
From Lemma 10,
$$ {|\Theta_{i}' (\xi)| \over |\Theta_{i}' (\eta)|}\leq C_{3}^{n_{0}}$$
for all $\xi$ and $\eta$ in $D^{+}_{i}$, and from Lemma 9,
$$ {|\Pi_{i}' (\xi)| \over |\Pi_{i}' (\eta)|}\leq C_{2}$$
for all $\xi$ and $\eta$ in $D^{+}_{i}$.
Again, because $\Theta_{i}$ is a schlicht function on $\bold C_{I}$ and
contracts the hyperbolic distance $d$ on $\bold C_{I}$, we have
$\Theta_{i}(D^{+}_{i})
\subset D^{+}_{i}$. Therefore, from the chain rule,
$$ {|\Sigma_{i}' (\xi)| \over |\Sigma_{i}' (\eta)|}\leq
C_{5}= C_{2}\cdot C_{3}^{n_{0}}$$
for all $\xi$ and $\eta$ in $D^{+}_{i}$.
\hfill \bull

\vskip5pt
\proclaim Lemma 12. Suppose $G$ is in ${\cal E}_{0}$ and $D$ is a
closed simply connected convex domain with $I\subset D \subset {\bf
C}_{G}$.  Then for all $z=x+yi$ with $y>0$ in $D$
$$\sin \Big( \arg (G(z))\Big)\geq {\sin (\arg (z)) \over N_{0}}$$
where $N_{0} =\sup_{\xi , \eta \in D}
|G'(\xi)/G '(\eta)|$ measures the distortion of $G$ on $D$.

\vskip5pt
{\bf Proof.} Since $G$ maps ${\bf UH}$ into itself, it contracts
the hyperbolic metric $dz/y$ on ${\bf UH}$. Suppose $z=x+yi$
with $y>0$ and $G (z)=X+Yi$. Then $|G '(z)| y \leq Y$.
Therefore,
$$ \sin (\arg (z)) = {y \over |z|} \leq {Y\over |z|\cdot
|G '(z)|} = \sin \Big( \arg (G(z))\Big) \cdot {|G (z)|\over |z|
\cdot |G'(z)|}.$$ So
$$ {\sin (\arg (z)) \over N_{0}} \leq \sin \Big( \arg (G (z))\Big)$$
for all $z=x+yi$ with $y>0$ in $D$.
\hfill \bull

\vskip5pt
Now we complete the proof of Theorem 18 as follows.

\vskip5pt
{\bf Proof of Theorem 18.}
For any $z_{0}$ in ${\bf UH}$, let $z_{i+1}= L_{i}(G_{i}(z_{i}))$ for
$0\leq i\leq n$. Since $0\leq \arg (z_{1}) \leq \pi/\gamma $,
the smallest positive integer $i$ such that $z_{i}$ lies in $\Delta_{i}$
must be either bigger than zero or not exist.
If such a positive integer does not exist, then $0\leq \arg (z_{n+1}) \leq
\pi -\theta_{n} \leq \pi -\sigma$.
Now let us suppose that this smallest number exists and is $i_{0}+1$.
Then
$z_{i_{0}}$ is in $A_{i_{0}}\not= \emptyset$ which is a subset of $D^{+}_{i_{0}}$.
Since $0< \arg (z_{i_{0}}) \leq \pi -\theta_{i}$, Lemmas 11 and
12 assure us there is a constant angle $0< \theta \leq \sigma$ such that
$0< \arg (z_{n+1}) \leq \pi -\theta$.
\hfill \bull

\section {Feigenbaum-Like Quadratic-Like Maps}

In this section, we discuss Feigenbaum-like quadratic-like maps and 
prove the result of Sullivan
which says that a Feigenbaum-like quadratic-like map has the {\sl a priori} 
complex bounds and is unbranched.

Let us first recall some facts about infinitely renormalizable folding
mappings.
Let $(U,V,f)$ be a real quadratic-like map,
that is, $f(U\cap \bold R^{1}) \subseteq V\cap \bold R^{1}$.
Conjugating by a linear fraction
transformation, we may assume that $f(-1)=f(1) = -1$ and that $f|[-1,
1]$
is a folding map of $[-1,1]$ with a unique quadratic critical point $0$.
For example, $P_{t}(z)=t-(1+t)z^{2}$ for $0\leq t \leq 1$ is a real
quadratic-like map whenever it is
restricted to any domain bounded by an equipotential curve.
Furthermore, suppose $(U,V,f)$ is infinitely renormalizable. Then
the filled-in Julia set $K_{f}$ is connected. Let $\beta_{0} =-1$
and let $\alpha_{0}$ be the fixed point of $f$ in $(-1, 1)$. Then
$\beta_{0}$ is
the non-separate fixed point, $\alpha_{0}$ is the separate fixed point,
and $K_{f}\cap \bold R^{1} =[-1,1]$.
The mapping $f_{0}=f|[-1,1]$ is
an infinitely $(n_{1}, n_{2}, \ldots , n_{k}, \ldots )$-renormalizable
folding mapping where $\{ n_{k}\}_{k=1}^{\infty}$ is the
maximum sequence
of integers such that $f_{0}$ is renormalizable about $m_{k}$ for
$m_{k}=\prod_{i=1}^{k}n_{i}$.
Let $I_{k}=[-a_{k}, a_{k}]$ be the maximal interval containing $0$
(set $m_{0}=0$ and $I_{0}=[-1,1]$) such that

\vskip5pt
{\bf (a)} $f_{0}^{\circ m_{k}}$ is monotone when restricted to
$[-a_{k}, 0]$ and to $[0, a_{k}]$,

\vskip5pt
{\bf (b)} $f_{0}^{\circ m_{k}}(I_{k}) \subset I_{k}$,

\vskip5pt
{\bf (c)} $I_{k}$, $f_{0}(I_{k})$, $\ldots$, $f_{0}^{\circ
(m_{k}-1)}(I_{k})$ have pairwise disjoint interiors, and

\vskip3pt
{\bf (d)} $f_{0}^{\circ m_{k}}$ has exactly two fixed points $\beta_{k}$
and $\alpha_{k}$ in $I_{k}$ where $a_{k}=|\beta_{k}|$.

\vskip3pt
\noindent There is a domain $U_{k}\supseteq I_{k}$ such that
$(U_{k}, V_{k}, f^{\circ m_{k}})$ is a quadratic-like map with connected
filled-in Julia set $K_{k}$, where $K_{k}$ is the
$k^{th}$-renormalization of $K_{1}$. Then
$\beta_{k}$ is the non-separate fixed point and $\alpha_{k}$ is the
separate fixed point of this quadratic-like map and then $K_{k}\cap {\bf
R}^{1} =I_{k}$. An infinitely renormalizable real quadratic-like map is
called Feigenbaum-like if $n_{k}=2$ for all $k$.

Suppose that $c(i)=f^{\circ i}(0)$ is the $i^{th}$ critical value of $f$
and that $J_{k}(i)$ is
the interval bounded
by $c(i)$ and $c(m_{k}+i)$ for $k\geq 0$ and $0\leq i<m_{k}$. We note
that $J_{k}(0)=J_{k}(m_{k})$. Then
$f_{0}:J_{k}(0)\rightarrow J_{k}(1)$ is folding for all $k\geq 0$ and
$f_{0}:J_{k}(i)\rightarrow J_{k}(i+1)$ is
a homeomorphism for every $k\geq 1$ and $1\leq i < m_{k}$.
Let $\zeta_{k} =\{ J_{k}(i)\}_{0\leq i <m_{k}}$ for $k\geq 0$.
Let $I_{k}(i)=f_{0}^{\circ i}(I_{k})$, $0\leq i<
m_{k}$, and let
$\xi_{k} =\{ I_{k}(i)\}_{0\leq i <m_{k}}$. Note that $J_{k}(i) \subseteq I_{k}(i)$.
We use $LI_{k}(i)$ and
$RI_{k}(i)$
to denote the intervals in $\xi_{k}$
adjacent to $I_{k}(i)$ and on the left and right sides of $I_{k}(i)$,
respectively
(there is only $LI_{k}(1)$ or
$RI_{k}(2)$ in $\xi_{k}$).
Let $LI^{+}_{k}(i)$ be the smallest interval
containing $LI_{k}(i)$ and the left end-point of $I_{k}(i)$ and
let $RI^{+}_{k}(i)$ be the smallest interval
containing $RI_{k}(i)$ and the right end-point of $I_{k}(i)$,
for $i =0$ or $3\leq i < m_{k}$. Let $LI^{+}_{k}(2) =[-1,
c(2)]$ and
$RI_{k}^{+}(1)=[c(1), 1]$. Similarly, we can define
$LJ_{k}(i)$ and $RJ_{k}(i)$ and
$LJ_{k}^{+}(i)$ and $RJ_{k}^{+}(i)$
for $0\leq i< m_{k}$. The following theorem is due to
Sullivan.

\vskip5pt
\proclaim Theorem 19~[SU2]. There is a constant
$C>0$ such that
$$ \min \{ |LI_{k}^{+}(i)|, |RI_{k}^{+}(i)| \} \geq C\cdot |I_{k}(i)|,$$
and such that
$$ \min \{ |LJ_{k}^{+}(i)|, |RJ_{k}^{+}(i)| \} \geq C\cdot |J_{k}(i)|.$$
for all $k\geq 0$ and $0\leq i <m_{k}$.

\vskip5pt
Consider the slit domain $V_{0}=V\setminus [c(1), \infty)$. The map
$f|V_{0}$ has two inverse branches (see Fig. 5)
$$g_{0} : V_{0} \rightarrow U_{0,0}= U\cap \{z=x+yi\ \in {\bf
C}\; |\; x<0\}$$
and
$$g_{1}: V_{0}\rightarrow U_{0,1}= U\cap \{z=x+yi\ \in {\bf
C}\; |\; x>0\}.$$

\vskip5pt
\centerline{\psfig{figure=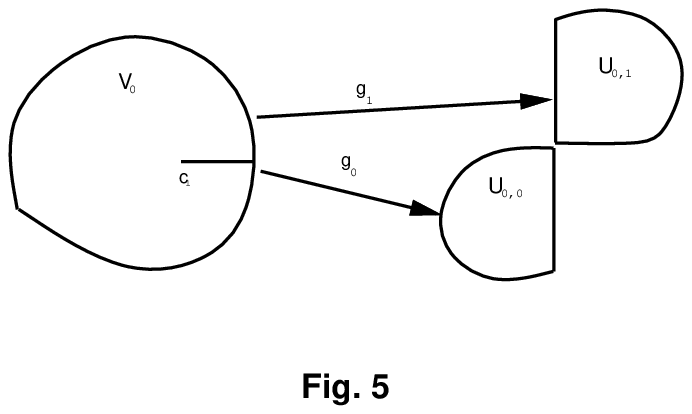}}
\vskip5pt

\noindent For each $k\geq 1$, let $g_{s(i)} : J_{k}(i+1) \rightarrow
J_{k}(i)$ for $1\leq i<m_{k}$ where $s(i) = 0$ if $J_{k}(i)$ is contained in
the negative half of the real line, and otherwise $s(1) =1$. Then
$${\cal A}_{k} =g_{s(1)} \circ g_{s(2)} \circ \cdots \circ g_{s(m_{k}-1)}:
J_{k}(m_{k})=J_{k}(0) \rightarrow J_{k}(1)$$
is a
homeomorphism and can be
extended homeomorphically to the maximum closed interval
$$T_{k}(m_{k})\supseteq LJ^{+}_{k}(m_{k})\cup J_{k}(m_{k})\cup
RJ^{+}_{k}(m_{k}).$$
Furthermore, ${\cal A}_{k}$ can be extended analytically to
$V_{k} = (V \setminus \bold R^{1}) \cup \TT_{k}(m_{k}).$
Let us continue to use ${\cal A}_{k}$ to
denote this extension. Let $U_{k}'= {\cal A}_{k}(V_{k})$, and let $U_{k}=f^{-1}(U_{k}')$ be the
pre-image of $U_{k}'$ under $f$. Since $f$ has no attractive and no
parabolic
periodic point, then $g_{0}\big( {\cal A}_{k} (T_{k}(m_{k})) \big)$ and
$g_{1}\big( {\cal A}_{k} (T_{k}(m_{k})) \big)$ are contained strictly in
$T_{k}(m_{k})$. Thus, $\overline{U}_{k} \subset V_{k}$ and
$$f^{\circ m_{k}}: U_{k}\rightarrow V_{k}$$
is a quadratic-like map (see Fig. 6).

\vskip5pt
\centerline{\psfig{figure=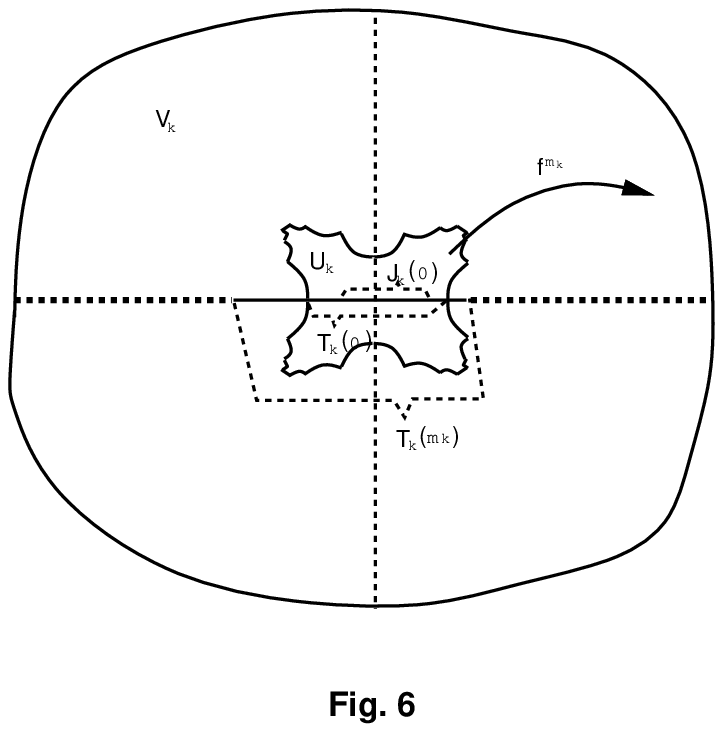}}
\vskip5pt

Let $T_{k}(i)=g_{s(i)}(T_{k}(i+1))$ for $i=m_{k}-1$, $m_{k}-2$,
$\ldots$, $1$. Let $T_{k}(0)= f^{-1}_{0} (T_{k}(1))$. Then
$$ J_{k}(0)=J_{k}(m_{k}) \subseteq I_{k}\subseteq T_{k}(0) \subset T_{k}(m_{k}).$$
The interval
$T_{k}(m_{k})$ is bounded by two critical values $c(q(k))$ and
$c(r(k))$
of $f$; one of them is a maximum value and the other is a minimum
value of $f^{\circ m_{k}}$. Suppose $T_{k}(1)=[d_{1}, e_{1}]$ where
$d_{1}< c(1)< e_{1}$. Let $d_{i}=f^{\circ i}(d_{1})$ and let $e_{i}
=f^{\circ i}(e_{1})$ for $1\leq i \leq m_{k}$. Note that $d_{m_{k}}$ and
$e_{m_{k}}$ are $c(r(k))$ and $c(q(k))$. We normalize $T_{k}(i)$ into
the
unit interval $[0,1]$: let $l_{i}: T_{k}(i) \rightarrow [0, 1]$ be the
linear map such that $l_{i}(d_{i}) =0$ and $l_{i}(e_{i})=1$ and let
$$\tilde{g}_{s(i)} = l_{i}\circ g_{s(i)}\circ l_{i+1}^{-1}$$
for $i=m_{k}-1$, $\ldots$, $2$, $1$.
The map $\tilde{g}_{s(i)}$ fixes $0$ and $1$ and
is a univalent function defined on
$$\bold C_{[0,1], V} = \Big( \big( \bold C\setminus \bold R^{1}\big)\cup
(0,1)\Big) \cap V.$$
The restriction
$\tilde{g}_{s(i)}|[0,1]$ is a homeomorphism of $[0,1]$.
Let
$${\cal L}_{k} =\tilde{g}_{s(1)}\circ \tilde{g}_{s(2)} \circ \cdots
\circ \tilde{g}_{s(m_{k}-1)}$$
for $k\geq 1$.

\vskip5pt
\proclaim Lemma 13. Suppose $(U,V,f)$ is a Feigenbaum-like quadratic-like
map. Then 
$$\{ {\cal L}_{k}\}_{k=3}^{\infty}$$ are uniform root-like maps.

\vskip5pt
{\bf Proof.} Since $f$ is a Feigenbaum-like map, we have
$m_{k}=2^{k}$ and $T_{k}(m_{k}) =J_{k-2}(0)$ for $3\leq k < \infty$.
For any $k\geq 3$, consider the homeomorphism
$f^{\circ (2^{k-1}-1)}: J_{k}(2^{k-1}+1) \rightarrow
J_{k}(2^{k})=J_{k}(0)$.
Its inverse $g_{s(2^{k-1}+1)} \circ \cdots \circ g_{s(2^{k}-1)} : J_{k}(2^{k})
\rightarrow J_{k}(2^{k-1}+1)$ can be extended to $T_{k-1}(2^{k-1})$, that
is, we can consider
$$g_{s(2^{k-1}+1)} \circ \cdots \circ g_{s(2^{k}-1)} : T_{k-1}(2^{k-1})
\rightarrow T_{k-1}(1).$$
Let
$$G_{0} = \tilde{g}_{s(2^{k-1}+1)} \circ \cdots \circ
\tilde{g}_{s(2^{k}-1)}$$
and let
$$L_{a_{0}} =\tilde{g}_{s(2^{k-1})}.$$
Then $G_{0}$ is a univalent map (or holomorphic embedding) of ${\bf
C}_{[0, 1], V}$ such that
$G_{0}(0)=0$ and $G_{0}(1)=1$ and such that $G_{0}|[-1,1]$ is a
homeomorphism of $[0,1]$. The map $L_{a_{0}}$
is a square root at $0$ (see Fig. 7).

Consider the homeomorphism
$f^{\circ (2^{k-2}-1)}: J_{k}(2^{k-2}+1) \rightarrow J_{k}(2^{k-1})$. Its inverse
$g_{s(2^{k-2}+1)} \circ \cdots \circ g_{s(2^{k-1}-1)} : J_{k}(2^{k-1})
\rightarrow J_{k}(2^{k-2}+1)$ can be extended to $T_{k-2}(2^{k-2})$, that
is, we can consider
$$G_{1}'= g_{s(2^{k-2}+1)} \circ \cdots \circ g_{s(2^{k-1}-1)} : T_{k-2}(2^{k-2})
\rightarrow T_{k-2}(1).$$
Let
$$G_{1} = \tilde{g}_{s(2^{k-2}+1)} \circ \cdots \circ
\tilde{g}_{s(2^{k-1}-1)}$$ and let
$$L_{a_{1}}=\tilde{g}_{s(2^{k-2})}.$$
The map $G_{1}$ is a univalent map of $\bold C_{[0, 1], V}$ and
$G_{1}|[-1,1]$ is a homeomorphism of $[0,1]$. The map $L_{a_{1}}$
is a square root at $a_{1}$.
The pre-image of $c(1)$ under $G_{1}'$
is $c(2^{k-2})$.
One of the end-points of $T_{k}(2^{k-1})$ is $0$; the other one is
between $J_{k-1}(0)$ and one of $LJ_{k-1}(0)$ or $RJ_{k-1}(0)$.
From Theorem 19, two components
of $T_{k-2}(2^{k-2})\setminus T_{k}(2^{k-1})$ have lengths greater than a
constant $C$ (obtained from Theorem 19)
times the length of $T_{k}(2^{k-1})$.
Thus
$(L_{a_{1}}, G_{1})$ is a
compatible pair and satisfies {\bf (i)} and {\bf (iii)} of Definition
4 (see Fig. 7).

Next we consider the homeomorphism
$f^{\circ (2^{k-3}-1)}: J_{k}(2^{k-3}+1) \rightarrow J_{k}(2^{k-2})$. Its
inverse $g_{s(2^{k-3}+1)} \circ \cdots \circ g_{s(2^{k-2}-1)} : J_{k}(2^{k-2})
\rightarrow J_{k}(2^{k-3}+1)$ can be extended to $T_{k-3}(2^{k-3})$;
that is, we can consider
$$ G_{2}' = g_{s(2^{k-3}+1)} \circ \cdots \circ g_{s(2^{k-2}-1)} : T_{k-3}(2^{k-3})
\rightarrow T_{k-3}(1).$$
Let
$$G_{2} = \tilde{g}_{s(2^{k-3}+1)} \circ \cdots \circ
\tilde{g}_{s(2^{k-2}-1)}$$
and let
$$L_{a_{2}}=\tilde{g}_{s(2^{k-3})}.$$
The map $G_{2}$ is a univalent map of $\bold C_{[0, 1], V}$ and
$G_{2}|[-1,1]$ is a homeomorphism of $[0,1]$. The map $L_{a_{2}}$
is a square root at $a_{2}$.
The preimage of $c(1)$ under $G_{2}'$ is $c(2^{k-3})$.
The interval $T_{k}(2^{k-2})$ is contained in $T_{k-1}(2^{k-1})$.
From Theorem 19, two components of $T_{k-3}(2^{k-3}) \setminus
T_{k}(2^{k-2})$
have lengths greater than a constant $C$ (obtained from Theorem 19)
times the length of $T_{k}(2^{k-2})$. Thus $(L_{a_{2}}, G_{2})$
is a compatible pair and satisfies {\bf (iii)}
of Definition 4 (see Fig. 7).

In general, for $3 < i \leq k-1$, consider the homeomorphism
$f^{\circ (2^{k-i}-1)}: J_{k}(2^{k-i}+1) \rightarrow J_{k}(2^{k-i+1})$.
Its inverse $g_{s(2^{k-i}+1)} \circ \cdots \circ g_{s(2^{k-i+1}-1)}
: J_{k}(2^{k-i+1}) \rightarrow J_{k}(2^{k-i}+1)$
can be extended to $T_{k-i}(2^{k-i})$; that is, we can consider
$$G_{i+1}' = g_{s(2^{k-i}+1)} \circ \cdots \circ g_{s(2^{k-i+1}-1)}
: T_{k-i}(2^{k-i}) \rightarrow T_{k-i}(1).$$
Let
$$G_{i+1} = \tilde{g}_{s(2^{k-i}+1)} \circ \cdots \circ
\tilde{g}_{s(2^{k-i+1}-1)}$$
and let
$$L_{a_{i+1}}=\tilde{g}_{s(2^{k-i})}.$$
The map $G_{i+1}$ is a univalent map of $\bold C_{[0, 1], V}$ and
$G_{i+1}|[0,1]$ is a homeomorphism. The map $L_{a_{i+1}}$
is a square root at $a_{i+1}$.
The preimage of $c(1)$ under $G_{i+1}'$ is $c(2^{k-i})$.
The interval $T_{k}(2^{k-i+1})$ is contained in
$T_{k-i+2}(2^{k-i+2})$.
From Theorem 19, two components
$T_{k-i}(2^{k-i}) \setminus T_{k}(2^{k-i+1})$
have lengths greater than a constant $C$
times the length of $T_{k}(2^{k-i+1})$.
Thus $(L_{a_{i+1}}, G_{i+1})$ is a
compatible pair and satisfies {\bf (iii)} of Definition 4.

One of the endpoints of $T_{k}(2^{k-i}+1)$ is to the left of
$c(2^{k-i}+1)$;
the other is in the interval in $\zeta_{k}$ which is adjacent to
$J_{k}({2^{k-i}+1})$.
The branch point of $g_{s(2^{k-i})}$ is always $c(1)$.
From
Theorem 19, we have a constant $\lambda >1$ such that $\{
a_{i}\}_{i=1}^{k}$ satisfies {\bf (ii)} of Definition 4 for $C=1$.

\vskip5pt
\centerline{\psfig{figure=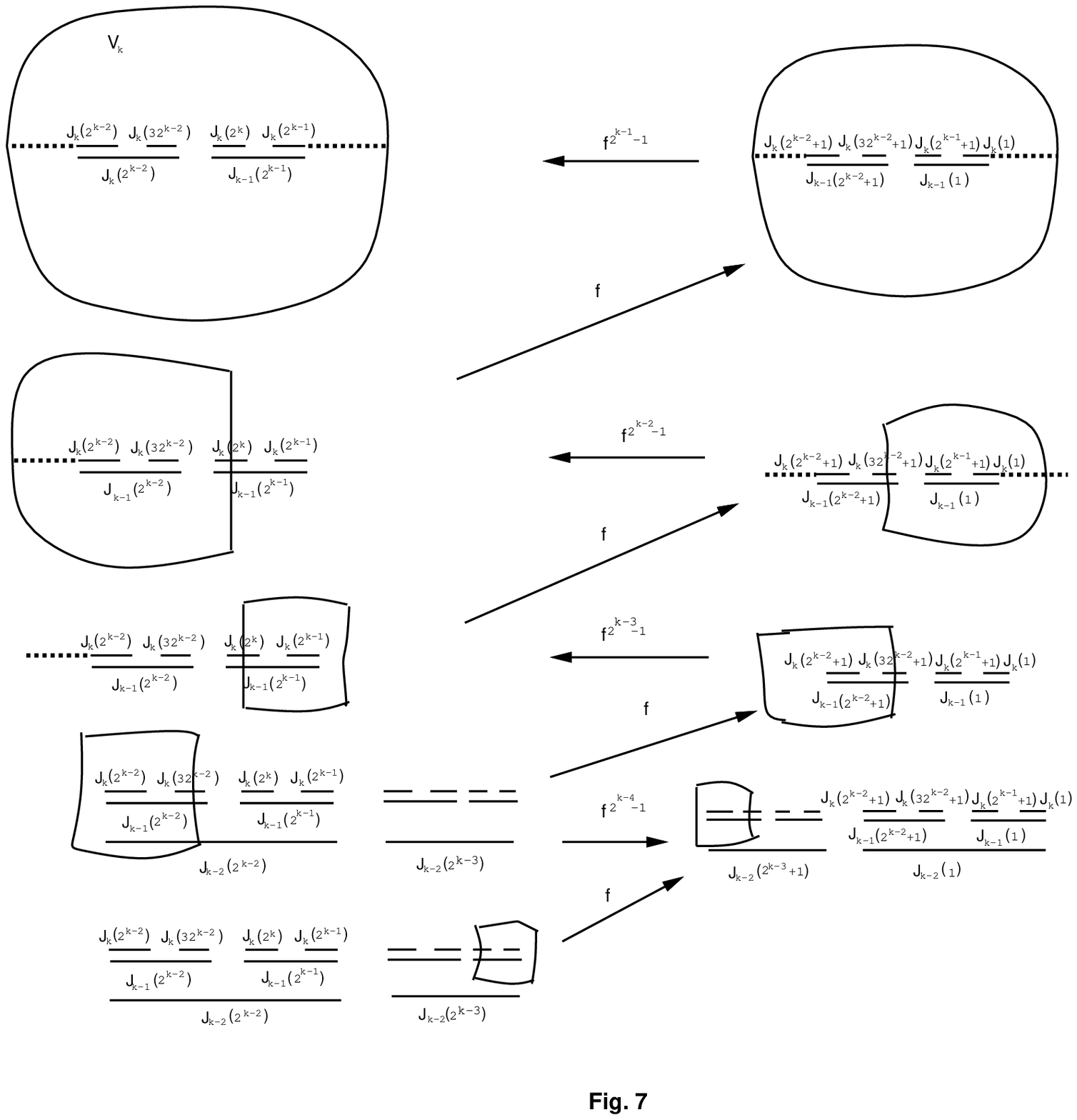}}
\vskip5pt

For all $k\geq 3$, we therefore decompose
$${\cal L}_{k}=L_{a_{k}}\circ G_{k}\circ L_{a_{k-1}}\circ G_{k-1}\circ \cdots \circ
L_{a_{1}}\circ G_{1}\circ L_{a_{0}}\circ G_{0}.$$
From the construction above,
${\cal L}_{k}$, $3\leq k< \infty$, is
uniform root-like map.
\hfill \bull

\vskip5pt
Let $(U,V,f)$ be a Feigenbaum-like quadratic-like map. Let
$m_{k}=2^{k}$.
Consider the renormalizations
$$f^{\circ m_{k}} : U_{k}\rightarrow V_{k}$$
for $1\leq k< \infty$ where $V_{k} = (V \setminus \bold R^{1}) \cup
\TT{k}(m_{k}).$
Let $U_{k}\cap \bold R^{1} =[-o_{k}, o_{k}]$. Let $w_{k, \theta}$ be the
the ray starting at $o_{k}$ with slop $\tan \theta$ for $0< \theta <\pi/2$.
Let $0\in R_{k, \theta}$ be the domain bounded by
$w_{k, \theta}$, $-w_{k, \theta}$, $\overline{w}_{k, \theta}$,
and $-\overline{w}_{k, \theta}$.

\vskip5pt
\proclaim Lemma 14. Suppose $(U,V,f)$ is a Feigenbaum-like
quadratic-like map. Then
there is a constant angle $\theta_{0} >0$ such that
$$ U_{k} \subseteq R_{k, \theta_{0}}$$
for all $k\geq 0$ (see Fig. 8).

\vskip5pt
\centerline{\psfig{figure=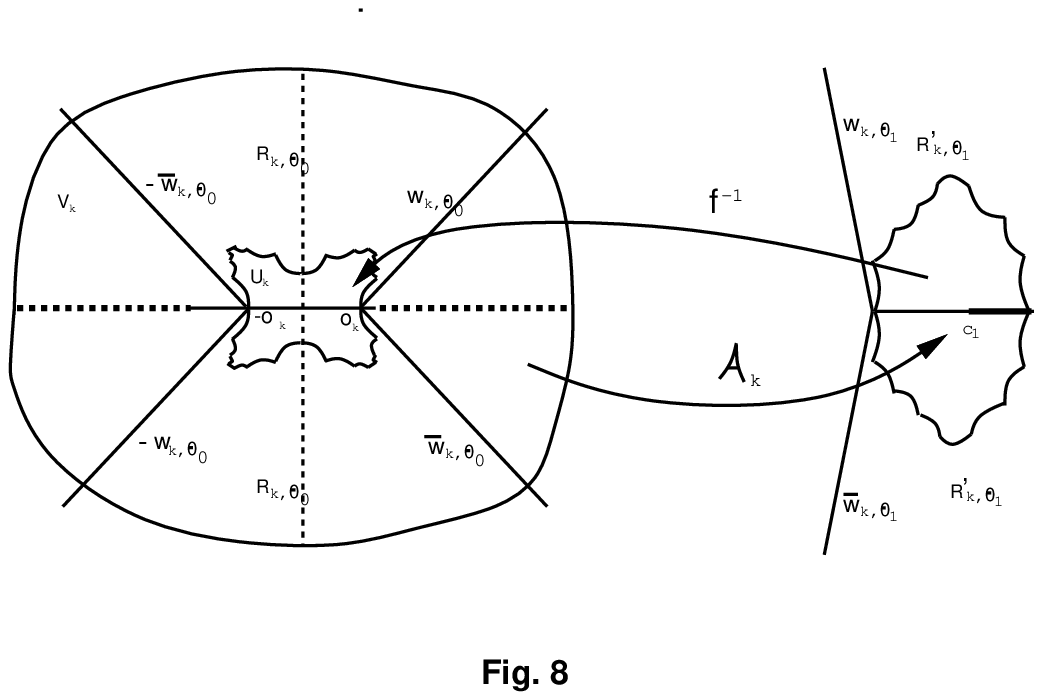}}
\vskip5pt

{\bf Proof.} Consider $f^{\circ (m_{k}-1)} : J_{k}(1) \rightarrow J_{k}(0)$.
Its inverse has the maximum extension
$${\cal A}_{k} = g_{s(1)} \circ g_{s(2)} \circ \cdots \circ g_{s(m_{k}-1)}:
V_{k} = (V \setminus \bold R^{1}) \cup \TT_{k}(m_{k}) \rightarrow U_{k}'.$$
Remember that
$${\cal L}_{k} =l_{1}\circ {\cal A}_{k}\circ l_{m_{k}}^{-1}$$
where $l_{m_{k}}$ and $l_{1}$ are the linear maps normalizing $T_{m_{k}}$
and $T_{k}(1)$, respectively, to $[0, 1]$.
Let
$$w'_{k, \theta} =\{ z\in \bold C \; | \; \arg (z-d_{1}) =\theta, \Im (z)
>0\}$$
be the ray starting at $d_{1}$ with angle $0\leq \theta \leq \pi $
where $U_{k}'\cap \bold R^{1} = T_{k}(1) = [d_{1}, e_{1}]$
with $d_{1} < c(1) < e_{1}$.
Let $R_{k, \theta}'$ be the sector containing $c(1)$ bounded by $w_{k,
\theta}'$ and $\overline{w'}_{k, \theta}$.
Applying Theorem 18, there is a constant angle $0< \theta_{1}\leq \pi$
such that
$U_{k}'$ is contained in a sector domain $R_{k, \theta_{1}}'$ (see Fig.
8). Since $f:
U\rightarrow V$ is a quadratic-like map, it is comparable with $z\mapsto
z^{2}$ near $0$.
So there is a constant angle $0< \theta_{0}< \pi /2$ depending on
$\theta_{1}$ such that $U_{k} =f^{-1} (U_{k}')$ is contained in $R_{k,
\theta_{0}}$ (see Fig. 8).
\hfill \bull

\vskip5pt
Take $I_{0}=(-1,1)$. Let $\bold C_{I_{0}} = (\bold C\setminus {\bf
R}^{1})\cup I_{0}$. Let $d=d_{H, \bold C_{I_{0}}}$ be the
hyperbolic distance on $\bold C_{I_{0}}$. Let
$\Omega_{r}=\{ z\in \bold C_{I_{0}}\; |\; d(z,
I_{0}) <r\}$ be a hyperbolic neighborhood.  From Lemma 6, $\Omega_{r}$
is the union of two disks $D^{+}_{\beta}$ and $D^{-}_{\beta}$ centered
at $c^{+}_{\beta}=i\cot \beta$ and $c^{-}_{\beta} =-i\cot\beta$ with
radii
$R^{+}_{\beta}=R^{-}_{\beta}= 1/\sin\beta$ where $\beta$ is the angle
between $\partial
D^{+}_{\beta}$ and the line $[1,\infty)$ at $1$ (see Fig. 9). Using
the law of cosines for the triangle $\Delta (c0z)$, for any
$z=re^{i\phi}$ in $\partial D^{+}_{\beta}$, $$r= \cot\beta  \sin \phi
+\sqrt{\csc^{2}\beta-\cot^{2}\beta \cos^{2}\phi}.$$

Let $q(z) =\sqrt{z}$ be
the square root from $\bold C\setminus \{ x<0\}$ to the right half-plane
${\bf RH}$. Let $\Pi_{\beta}^{+} =q(D^{+}_{\beta})$. For
$0<\beta <\pi/4$, consider $\Omega_{r'}=D^{+}_{2\beta}\cup
D^{-}_{2\beta}$.
Let $z_{0}=re^{i\tau}\neq 1$ be the intersection point
of $\partial D^{+}_{2\beta}$ and $\Pi^{+}_{\beta}$. Then $\tau$
is the unique non-zero solution of the equation
$$
\displaylines{\sqrt{\cot\beta  \sin (2\phi)
+\sqrt{\csc^{2}\beta-\cot^{2}\beta \cos^{2}(2\phi)}}\hfill\cr=\hfill{}
\cot(2\beta) \sin \phi +\sqrt{\csc^{2}(2\beta)-\cot^{2}(2\beta)
\cos^{2}\phi}.}
$$
Thus $\tau=\tau (\beta)$ is a strictly increasing function and
$\tau\rightarrow 0$ as $\beta \rightarrow 0$.
Let $z_{0}-1=\tilde{r} e^{i\theta}$. Then $\theta=\theta
(\tau)$ is a strictly increasing function and
$\theta\rightarrow 0$ as $\tau\rightarrow 0$.
Let
$\theta=\theta \circ \tau (\beta)$. It is a strictly
decreasing function and
$\theta \rightarrow 0$ as $\beta \rightarrow 0$. Let $\beta =\beta
(\theta)$ be its inverse function.

\vskip5pt
\centerline{\psfig{figure=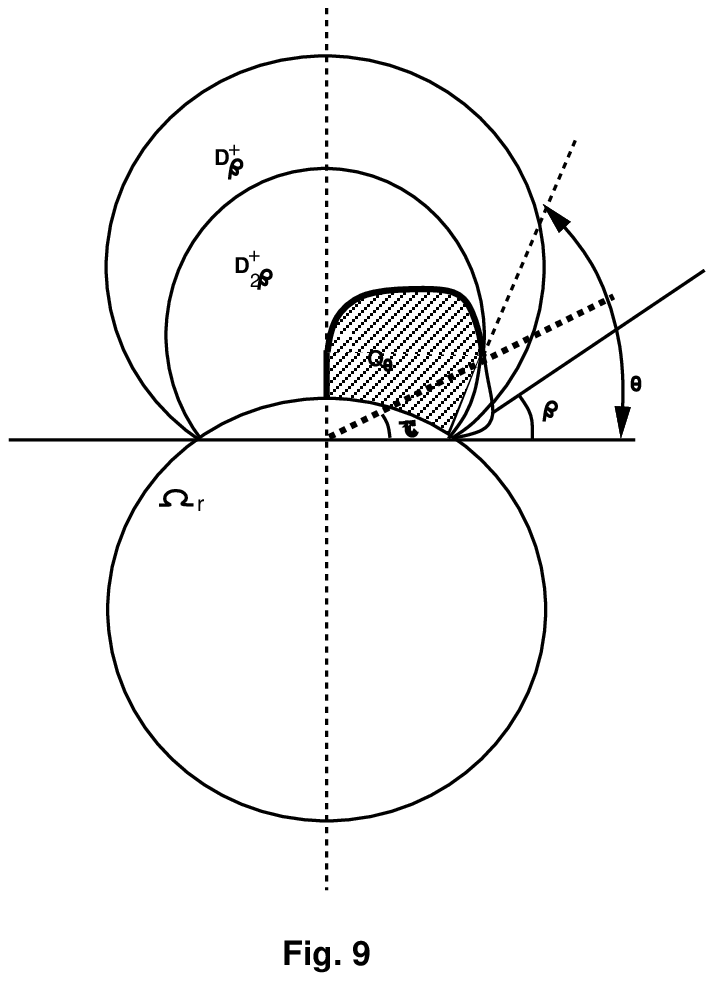}}
\vskip5pt

For any $0< \theta_{0}<\pi/2$, let $0< \beta_{0}=\beta
(\theta_{0})< \pi/2$ and let
$${\cal Q}_{\theta_{0}} =\{ z\in \Pi_{\beta_{0}}
\;|\; \theta_{0}\leq \arg (z-1) \leq {\pi\over 2}\}.$$
Then ${\cal Q}_{\theta_{0}}\subseteq D^{+}_{2\beta_{0}}$. Let
$\overline{\cal Q}_{\theta_{0}}=\{ \overline{z} \; |\; z\in
{\cal Q}_{\theta_{0}}\}$
and let ${\cal S}_{\theta_{0}} ={\cal Q}_{\theta_{0}} \cup
(-{\cal Q}_{\theta_{0}})\cup
\overline{\cal Q}_{\theta_{0}} \cup (-\overline{\cal
Q}_{\theta_{0}})$. For a
number $0< \nu_{0} <1$, let
$$\nu_{0} \cdot {\cal S}_{\theta_{0}} =\{ w=\nu_{0}\cdot
z \;| \; z\in {\cal S}_{\theta_{0}}\}.$$
Let $A_{\theta_{0}} =\Omega_{r}
\setminus (\nu_{0} {\cal S}_{\theta_{0}})$. From the above calculation,
we have

\vskip5pt
\proclaim Lemma 15. There is a constant $C=C(\theta_{0}, \nu_{0} )>0$
depending
on $\theta_{0}$ and on $\nu_{0}$ such that the modulus
$\hbox{mod} (A_{\theta_{0}})$ of the annulus $A_{\theta_{0}}$ is greater than $C$.

\vskip5pt
{\bf Proof.}  Let
$a=\hbox{diam}(\nu_{0} {\cal S}_{\theta_{0}})$ be the diameter of
$\nu_{0} {\cal S}_{\theta_{0}}$ and let
$b=d(\partial (\nu_{0} {\cal S}_{\theta_{0}}), \partial \Omega_{r})$
be the distance between
$\partial (\nu_{0} {\cal S}_{\theta_{0}})$ and $\partial \Omega_{r}$.
Then
$a/b$ is
bounded from above by a constant
depending only on $\theta_{0}$ and on $\nu_{0}$. This implies the lemma
(by using Gr\"otzsch argument (refer to~[AH1])).
\hfill \bull

\vskip5pt
Suppose $0< C_{0}<1$ is a constant. If we use $q_{a}(z)
=\sqrt{z-a}/\sqrt{1-a}$ for $|a|<C_{0}$ to replace $q(z)$
in the above calculation, then Lemma 15 has a generalized version.

\vskip5pt
\proclaim Lemma 15$^{\prime}$. There is a constant $C=C(\theta_{0}, \nu_{0}, C_{0}
)>0$ depending
on $\theta_{0}$, $\nu_{0}$, and $C_{0}$ such that the modulus
$\hbox{mod} (A_{\theta_{0}})$ of the annulus $A_{\theta_{0}}$ is greater
than $C$.

\vskip5pt
\proclaim Theorem 20~[SU2]. A Feigenbaum-like quadratic-like
map $(U, V, f)$ has the {\sl a priori} complex bounds and is unbranched.

\vskip5pt
{\bf Proof.} We use the same notation as in the previous lemmas.
Suppose $(U,V,f)$ is a Feigenbaum-like quadratic-like map. Let $m_{k}=2^{k}$.
Consider the renormalizations
$$f^{\circ m_{k}} : U_{k}\rightarrow V_{k}$$
for $1\leq k< \infty$ where
$V_{k} = (V \setminus \bold R^{1}) \cup \TT_{k}(m_{k})$ and
where $U_{k}\cap \bold R^{1} =[-o_{k}, o_{k}]$. From Lemma 14, there
is a constant angle $0< \theta_{0}< \pi/2$ such that
$$ U_{k} \subseteq R_{k, \theta_{0}}$$
for $k>0$. Let $\beta_{0}=\beta(\theta_{0}).$

Consider $f^{\circ (m_{k}-1)} : J_{k}(1) \rightarrow J_{k}(0)$.
Its inverse has the maximum extension
$${\cal A}_{k} = g_{s(1)} \circ g_{s(2)} \circ \cdots \circ g_{s(m_{k}-1)}:
V_{k} \rightarrow U_{k}'$$
where
$U_{k}'\cap \bold R^{1} = T_{k}(1) = [d_{1}, e_{1}]$
with $d_{1} < c(1) < e_{1}$. From Theorem 19 (and the bounded geometry
property of the attractive Cantor set), there is a constant $C_{0}>0$
such that $C_{0}^{-1} \leq |c(1)-d_{1}|/|e_{1}-c(1)| <C_{0}$ for all
$k>0$.
We normalize $\TT_{k}(m_{k})=(d_{m_{k}}, e_{m_{k}})$ to $(-1,1)$ by
the linear map $s_{1}$ such that $s_{1}(d_{m_{k}})=1$ and
$s_{1}(e_{m_{k}})=-1$. We normalize $(0, o_{k})$ to $(0, 1)$ by the
linear map $s_{2}$ such that $s_{2}(0)=0$ and $s_{2}(o_{k})=1$. We
normalize $\TT_{k}(1) =(d_{1}, e_{1})$ to $(-1, 1)$ by the linear
map $s_{3}$ such that $s_{3}(e_{1}) =-1$ and
$s_{3}(d_{1}) =1$. Let $a=s_{3}(c(1))$. Then there is a constant we
still denote it as
$0< C_{0}<1$ such that $|a| < C_{0}$ for all $k>0$.
There is an integer $n_{0}>0$ such that
for any $k>n_{0}$, $\Omega_{r}=D^{+}_{\beta_{0}}\cup D^{-}_{\beta_{0}}$ is
contained in $s_{1}(V)$. Let ${\cal B}_{k} =s_{3}\circ {\cal
A}_{k}\circ s_{1}^{-1}$ and let $q_{a}(z) =s_{2}\circ g_{s_{m_{k}}}\circ
s_{3}^{-1}$. Then $q_{a}$ is comparable with $\sqrt{z-a}/\sqrt{1-a}$.
Since ${\cal B}_{k}$ contracts the hyperbolic distance $d_{H,{\bf
C}_{I_{0}}}$, then ${\cal B}_{k}
(\Omega_{r})\subseteq \Omega_{r}$. From Lemma 14,
$$X_{k}'=q_{a}({\cal B}_{k} (\Omega_{r})) \subseteq {\cal
S}_{\theta_{0}} ={\cal Q}_{\theta_{0}} \cup (-{\cal Q}_{\theta_{0}})\cup
\overline{\cal Q}_{\theta_{0}} \cup (-\overline{\cal
Q}_{\theta_{0}}).$$
Let $X''_{k} =s_{1}\circ s_{3}^{-1}(X_{k}')$. Then $X_{k}''\subseteq \nu_{0}
S_{\theta_{0}}$ for all $k>0$ where $\nu_{0}>0$ is a constant obtained from
Lemma 4.2.
Let $Y_{k}''=\Omega_{r}$. Then
from Lemma 15$^{\prime}$,
$$\hbox{mod} (Y_{k}''\setminus X_{k}'') >C$$
for all $k>n_{0}$ where $C>0$ is a constant.

Now let $X_{k} =s_{1}^{-1} (X_{k}'')$ and
let $Y_{k}=s_{1}^{-1}(Y_{k}'')$.
Then
$$f^{\circ m_{k}} : X_{k}\rightarrow Y_{k}$$
is quadratic-like map and $\hbox{mod} (Y_{k}\setminus X_{k}) >C$ for
all $k>n_{0}$. This means that $(U, V, f)$ has the {\sl a
priori} complex bounds.

Let $W_{k} =Y_{k} \setminus (LJ_{k}(0)\cup RJ_{k}(0))$. Applying Theorem
19 and the above argument, there is a constant $C'>0$ such that
$\hbox{mod}(W_{k}\setminus
K_{k}) >C'$ for all $k>n_{0}$ where $K_{k}$ is the filled-in Julia set
of $f^{\circ m_{k}} : X_{k}\rightarrow Y_{k}$. But $W_{k}\cap CO =\{
c(jm_{k})\}_{j=0}^{\infty}$. So $(U, V, f)$ is unbranched. This
completes the proof.
\hfill \bull

\vskip5pt
Theorem 20 and Theorem 16 now give us that

\vskip5pt
\proclaim Corollary 1. The filled-in Julia set
$K_{f}$ of a Feigenbaum-like quadratic-like map $(U, V, f)$
is locally connected.

\vskip5pt
Sullivan~[SU2] (see also~[MEV]) also proved that any real infinitely
renormalizable quadratic-like map $f:U\rightarrow V$ of bounded type
has the {\sl a priori} complex bounds.
Thus it is unbranched and its filled-in Julia set is locally
connected from Theorem 16. 

\section {The Local Connectivity of Certain Infinitely Renormalizable 
Quadratic Julia Sets}

We prove in this section the following result.

\vskip5pt
\proclaim Theorem 21. There is a subset
$\tilde{\Upsilon}$ in
${\cal M}$ such that {\bf (1)} $\tilde{\Upsilon}$ is dense in the
boundary $\partial {\cal M}$ of the Mandelbrot set ${\cal M}$,
${\bf (2)}$ for every $c$ in $\tilde{\Upsilon}$, $P_{c}$ is unbranched,
infinitely renormalizable and has the {\sl a priori} complex bounds.

\vskip5pt
From Theorem 16, we have

\vskip5pt
\proclaim Corollary 2. The filled-in Julia set $K_{c}$ of
$P_{c}$ is locally connected for every $c$ in $\tilde{\Upsilon}$.

\vskip5pt
{\bf Proof of Theorem 21.}
Suppose $c_{0}$ is a Misiurewicz point in ${\cal M}$.
Then there is an integer $m>1$ such that $p=P_{c_{0}}^{\circ m}(0)$ is a
repelling periodic point of period $k\geq 1$. Let $\alpha$ be the separate
fixed point of $P_{c_{0}}$. Without loss of generality, we assume that
$P_{c_{0}}$ is non-renormalizable. (If $P_{c_{0}}$ is renormalizable,
it must be finitely renormalizable. We would then take $\alpha$ as
the separate
fixed point of the last renormalization of $P_{c_{0}}$ (see \S 4)).
Let $\Gamma$ be
the union of a cycle of external rays landing at $\alpha$. Let $\gamma$
be a fixed equipotential curve of $P_{c_{0}}$.
Using $\Gamma$ and
$\gamma$, we construct the two-dimensional Yoccoz puzzle as follows
(see \S 3).
Let $C_{-1}$ be the domain bounded by $\gamma$. The set $\Gamma$ cuts
$C_{-1}$ into a finite number of closed domains. Let $\eta_{0}$ denote the
set of these domains. Let $\eta_{n} =P_{c_{0}}^{-n} (\eta_{0})$.
Let $C_{n}$ be the member of $\eta_{n}$ containing $0$.

Let
$$ p\in \cdots \subseteq D_{n}(p) \subseteq D_{n-1}(p) \subseteq \cdots
\subseteq D_{1}(p) \subseteq D_{0}(p)$$
be a $p$-end, where $D_{n}(p) \in \eta_{n}$.
Let
$$ c_{0}\in \cdots \subseteq E_{n}(c_{0}) \subseteq E_{n-1}(c_{0})
\subseteq \cdots \subseteq E_{1}(c_{0}) \subseteq E_{0}(c_{0})$$
be a $c_{0}$-end, where $E_{n}(c_{0}) \in \eta_{n}$.
We have $P_{c_{0}}^{\circ (m-1)}(E_{n+m-1}(c_{0})) =D_{n}(p)$.

Since the diameter $\hbox{diam}(D_{n}(p))$ tends to zero
as $n\rightarrow \infty$ and since $p$ is a repelling periodic point,
we can find an integer $l\geq m$ such that $|(P_{c_{0}}^{\circ
k})'(x)|\geq
\lambda >1$ for all $x\in D_{l}(p)$ and such that
$$P_{c_{0}}^{\circ (m-1)}: E_{l+m-1}(c_{0}) \rightarrow D_{l}(p)$$
is a homeomorphism.
Let $q>0$ be the smallest integer such that $P_{c_{0}}^{\circ
q}(D_{l}(p))$ contains $0$, i.e., it is $C_{r_{0}}$ in $\eta_{r_{0}}$
where $r_{0}\geq 0$. Then
$$ f=P_{c_{0}}^{\circ q}: D_{l}(p) \rightarrow C_{r_{0}}$$
is a homeomorphism. Let $r>r_{0}$ be an integer such that
$B_{0}=f^{-1}(C_{r})\subset D_{l}(p)$ does not contain $p$ where $C_{r}$ is
the member of $\eta_{r}$ containing $0$.
Then
$$ P_{c_{0}}^{\circ q}: B_{0} \rightarrow C_{r}$$
is a homeomorphism.
The domain $B_{0}$ is a member of $\eta_{r+q}$.
Let $B_{n}\subseteq D_{l+nk}$ be the pre-image of $B_{0}$ under
$P_{c_{0}}^{\circ nk}|D_{l+nk}(p)$ for $n\geq 1$.
The domain $B_{n}$ is a member of $\eta_{r+q+nk}$ and
$$P_{c_{0}}^{\circ (q+nk)} : B_{n} \rightarrow C_{r}$$
is a homeomorphism.

From the structural stability theorem (see~[PRZ,SHU]), the points
$\alpha$ and $p$ and the sets $\Gamma$, $C_{r}$, $D_{n}$ and $B_{n}$ for
$n\geq 1$ are all preserved by a small perturbation $c$ of $c_{0}$
(refer to~[JI4]).
Therefore they can be constructed for $P_{c}$ as long as $c$ near $c_{0}$.
Let $U_{0}$ be a neighborhood about $c_{0}$
such that the corresponding points
$\alpha(c)$ and $p(c)$ and the corresponding sets $\Gamma (c)$,
$C_{r}(c)$, $D_{n}(c)$,
and $B_{n}(c)$ for $n\geq 1$ are all preserved for $c\in U_{0}$.
Moreover, as $n$ goes to infinity, the diameter $\hbox{diam}(B_{n}(c))$
tends to zero and the set $B_{n}(c)$ approaches to $p(c)$ uniformly on $U_{0}$.
Let
$$W_{n}=W_{n}(c_{0})=\{ c \in \bold C \; |\; P_{c}^{m} (0) \in
B_{n}(c)\}.$$
Then from the result in~[JI4], $W_{n}\subseteq U_{0}$ for $n$ large
enough.

For any $c\in W_{n}$, Let $R_{n}(c)$ be the preimage of $B_{n}(c)$ under the
map $P_{c_{0}}^{\circ (m-1)}: E_{l+m-1}(c) \rightarrow D_{l}(p,c)$ and let
$C_{r+q+nk+m}(c)=P_{c}^{-1}(R_{n}(c))$. The domain $C_{r+q+nk+m}(c)$ is
the member containing $0$ in $\eta_{r+q+m+nk}$ and
$$F_{n,c}=P_{c}^{\circ (q+nk+m)}:
X_{n}(c)=\CC_{r+q+nk+m}(c)\rightarrow Y_{n}(c) = \CC_{r}(c)$$
is a quadratic-like map. Let $A_{n}(c)= \CC_{r}(c)\setminus
C_{r+q+nk+m}(c)$. Since the diameter
$\hbox{diam}(C_{r+q+nk+m})(c))$ tends to zero as $n$ goes to
infinity uniformly in $U_{0}$, there
is an integer $N_{0}>0$ such that
$$\hbox{mod} (A_{n}(c)) \geq 1$$
for all $n\geq N_{0}$ and $c\in W_{n}$.
Since
$$ \{ F_{n,c}: X_{n}\rightarrow Y_{n} \; |\; c\in W_{n}\}$$
is a full family of quadratic-like maps, $W_{n}$ contains
a copy ${\cal M}_{n}={\cal M}_{n}(c_{0})$ of the Mandelbrot set ${\cal
M}$ (see~[DH3]). For any $c\in {\cal M}_{n}$, $P_{c}$ is once
renormalizable and
$$CO(c) \cap C_{r+q+nk+m}(c) =\{ c(j(q+nk+m))\}_{j=0}^{\infty}$$
where $CO(c)=\{ c(i)=P_{c}^{\circ i}(0)\}_{i=0}^{\infty}$.
Let
$$\tilde{\Upsilon}_{1}(c_{0}) =\cup_{n=N_{0}}^{\infty} {\cal M}_{n}.$$

We use the induction to complete the construction of the subset
$\tilde{\Upsilon}(c_{0})$ around $c_{0}$. Suppose we have constructed
$W_{w}$
where $w=i_{0}i_{1}\ldots i_{k-1}$ and $i_{0}\geq N_{0}$, $i_{1}\geq
N_{i_{1}}$, $\ldots$, $i_{k-1}\geq N_{i_{0}i_{1}\ldots i_{k-2}}$.
There is a parameter $c_{w} \in {\cal M}_{w}$ such
that $$F_{w}=F_{w, c_{w}}: X_{w}=X_{w}(c_{w})
\rightarrow Y_{w}=Y_{w}(c_{w})$$
is hybrid equivalent (see \S 1) to $P(z)=z^{2}-2$.
For $F_{w}$, let $\beta_{w}$ and $\alpha_{w}$ be its
non-separate and separate fixed points. Let
$\tilde{\beta}_{w}$ be another preimage of $\beta_{w}$ under
$F_{w}$. Let $\Gamma_{w}$ be the external rays of
$P_{c_{w}}$ landing at $\alpha_{w}$. Let $Y_{w0}$ be the
domain containing $0$ and bounded by $\partial X_{w}$ and
$F_{w}^{-1}(\Gamma_{w})$. Let $\tilde{\beta}_{w} \in
E_{w0}$ and $\beta_{w} \in E_{w1}$ be the components
of the closure of $X_{w}\setminus Y_{w0}$. Let $G_{w0}$ and
$G_{w1}$ be the inverses of $F_{w}|E_{w0}$ and
$F_{w}|E_{w1}$. Let
$$D_{wn}=G_{w1}^{\circ n}(D_{w0})$$
and
$$ B_{wn} = G_{w0}(D_{w(n-1)})$$
for $n\geq 1$.
From the structural stability theorem (see~[PRZ,SHU]),
the points $\beta_{w}$ and $\alpha_{w}$ and the sets $\Gamma_{w}$
are all preserved by a small perturbation $c$ of $c_{w}$. Therefore
we can find a small neighborhood $U_{w}$ about $c_{w}$ such that
the corresponding domains $D_{wn}(c)$ and $B_{wn}(c)$ can be constructed
for $P_{c}$, $c\in U_{w}$ (refer to~[JI4]). Let
$$W_{wn} =\{ c\in \bold C \; |\; F_{w,c}(0) \in B_{wn}(c)
\}.$$
The diameter $\hbox{diam}(B_{wn}(c))\rightarrow 0$
as $n\rightarrow \infty$ uniformly on $U_{w}$.
From the result in~[JI4], $W_{wn}\subseteq
U_{w}$ for $n$ large.

For each $c$ in $W_{wn}$, $n\geq N_{w}$,
let $X_{wn}(c) =F^{-1}_{w,c}(\BB_{wn}(c))$ and $Y_{wn}(c) =\YY_{w0}(c)$.
Then
$$F_{wn,c}=F^{\circ (n+1)}_{w,c} : X_{wn}(c) \rightarrow
Y_{wn}(c)$$
is a quadratic-like map.
Let
$$A_{wn}(c)= X_{wn}(c)\setminus \overline{Y}_{wn}(c).$$
Since the diameter $\hbox{diam}(Y_{wn}(c))$ tends to zero as $n$
goes to infinity uniformly in $U_{w}$. There
is an integer $N_{w}>0$ such that $$\hbox{mod} (A_{wn}(c)) \geq 1$$
for all $n\geq N_{w}$ and $c\in W_{wn}$.
Since
$$ \{ F_{wn,c}: X_{wn}\rightarrow Y_{wn} \; |\; c\in W_{wn}\}$$
is a full family of quadratic-like maps, $W_{wn}$ contains
a copy ${\cal M}_{wn}={\cal M}_{wn}(c_{0})$ of the Mandelbrot set
${\cal M}$ (see~[DH3]).
For any $c\in {\cal M}_{wn}$, $P_{c}$ is $k$-{\sl times} renormalizable
and $A_{wn}(c)$ contains no critical values of $P_{c}$.
Let
$$\tilde{\Upsilon}_{k}(c_{0})=\cup_{w}\cup_{n=N_{w}}^{\infty} {\cal
M}_{wn}$$ where $w$ runs over all sequences of integers of length $k$
in the induction.

We have thus constructed a subset
$\tilde{\Upsilon} (c_{0}) =\cap_{k=1}^{\infty}
\tilde{\Upsilon}_{k}(c_{0})$ such that every $c\in
\tilde{\Upsilon} (c_{0})$ is infinitely renormalizable and such
that $c_{0}$ is a limit point of $\tilde{\Upsilon} (c_{0})$.
From the above construction, for every $c\in \tilde{\Upsilon} (c_{0})$,
$P_{c}$ is unbranced and has the {\sl a priori} complex bounds.

Let $\tilde{\Upsilon} =\cup_{c_{0}}\tilde{\Upsilon} (c_{0})$ where
$c_{0}$ runs over all
Misiurewicz points in ${\cal M}$. Then for every $c\in
\tilde{\Upsilon}$, $P_{c}$ is
unbranched infinitely renormalizable and has the {\sl a priori} complex
bounds. Since the
set of Misiurewicz points is dense in $\partial {\cal M}$ (see~[CAG]),
the set $\tilde{\Upsilon}$ is dense in $\partial {\cal M}$.
It completes the proof of the theorem.
\hfill \bull

\vskip5pt
{\bf Remark 3.} Douady (see~[MI2]) constructed an example of an
infinitely
renormalizable quadratic polynomial whose filled-in Julia set is not locally
connected.

\end{document}